\documentclass[10pt]{article}
\expandafter\let\csname equation*\endcsname\relax
\expandafter\let\csname endequation*\endcsname\relax
\usepackage{bm}
\usepackage{enumerate}
\usepackage[top=1in, bottom=1in, left=1.2in, right=1.2in]{geometry}
\usepackage{mathrsfs,color}
\usepackage{amsfonts}\usepackage{multirow}
\usepackage{amssymb,dsfont}
\usepackage{commath}
\usepackage{caption,comment}
\usepackage{blkarray}
\usepackage{amsmath, amscd}
\usepackage{tikz-cd}
\usepackage{float}
\usepackage{amssymb,amsthm}
\usepackage[toc,page]{appendix}
\usepackage{graphicx,afterpage}
\usepackage{epstopdf}
\usepackage{chngcntr}
\usepackage{apptools}
\AtAppendix{\counterwithin{theorem}{section}}

\RequirePackage{doi}
\usepackage{hyperref}

\hypersetup{CJKbookmarks,%
	bookmarksnumbered,%
	colorlinks,%
	linkcolor=black,%
	citecolor=black,%
	plainpages=false,%
	pdfstartview=FitH}
\numberwithin{equation}{section}
\numberwithin{figure}{section}

\newcommand\tabcaption{\def\@captype{table}\caption}

\newtheorem{thm}{Theorem}[section]

\newtheorem{aspt}[thm]{Assumption}

\newtheorem{theorem}{Theorem}[section]
\newtheorem{lemma}[theorem]{Lemma}

\allowdisplaybreaks

 \usepackage{array}
 \newcolumntype{M}[1]{>{\centering\arraybackslash}m{#1}}
 \newcolumntype{N}{@{}m{0pt}@{}}
\usepackage[toc,page]{appendix}

\usepackage{authblk}

\date{}

\begin{document}
	
	\title{Analysis of the feedback particle filter with diffusion map based approximation of the gain}
	
%	\affil{First Institution/Department, Affiliation, City, Country }{FIRSTAFF}
%	\affil{Second Institution/Department, Affiliation, City, Country }{SECONDAFF}
%	
%	\author{A.B. Firstauthor}{FIRSTAFF}
%	\author{C. Coauthor}{SECONDAFF}
%	\author{D.E. Othercoauthor}{FIRSTAFF}

	\author[1]{Sahani Pathiraja}
	\author[2,3]{Wilhelm Stannat}
	\affil[1]{ \small{Institute of Mathematics, University of Potsdam, Germany}}
	\affil[2]{ \small{Institute of Mathematics, TU Berlin, Germany}}
	\affil[3]{ \small{Bernstein Center for Computational Neuroscience, Berlin, Germany}}
	
	\maketitle
	%\tableofcontents

	%%%%%%%%%%%%%%%%%%%%%%%%%%%%%%%%%%%%%%%%%%%%%%%%%%%%%%%%%%%%%%%%%%%%%%%%%%	
		\begin{abstract} 
		
			Control-type particle filters have been receiving increasing attention over the last decade as a means of obtaining sample based approximations to the sequential Bayesian filtering problem in the nonlinear setting.  Here we analyse one such type, namely the feedback particle filter and a recently proposed approximation of the associated gain function based on diffusion maps.  The key purpose is to provide analytic insights on the form of the approximate gain, which are of interest in their own right.  These are then used to establish a roadmap to obtaining well-posedness and convergence of the finite $N$ system to its mean field limit.  A number of possible future research directions are also discussed.   
		
	\end{abstract}
	%%%%%%%%%%%%%%%%%%%%%%%%%%%%%%%%%%%%%%%%%%%%%%%%%%%%%%%%%%%%%%%%%%%%%%%%%%%%%	
	\noindent
	{\bf Keywords.}  nonlinear filtering, diffusion maps, Poisson equation, propagation of chaos, stochastic analysis \\
	\noindent {\bf AMS(MOS) subject classifications.}  	60H10, 35Q93, 35J05, 60J60, 60J27, 62M05

	\vspace{0.5cm}

	\section{Introduction}
	The filtering problem is concerned with the following system
	\begin{align*}
		dS_t &= \mathcal{M}(S_t)dt + dV_t, \\
		dZ_t &= h(S_t)dt + dW_t,
	\end{align*}
	where $S_t \in \mathbb{R}^d$ and $Z_t \in \mathbb{R}$ with $\mathcal{M}: \mathbb{R}^d \rightarrow \mathbb{R}^d$ and $h: \mathbb{R}^d \rightarrow \mathbb{R}$.  $V_t$ and $W_t$ are independent Wiener processes and we define the filtration $\mathcal{Z}_t:= \sigma(Z_s: s \leq t)$.  The goal of the filtering problem is to estimate the conditional distribution $P\{S_t \in A | \mathcal{Z}_t\}$ (also referred to as the filtering distribution or posterior distribution in Bayesian terminology).  The Kushner-Stratonovich SPDE describes the time evolution of the filtering distribution.  It is well-known that there is no closed-form exact description of the conditional distribution for general $\mathcal{M}, h$, with the exception of the linear case with Gaussian initial condition, for which the Kalman-Bucy filter provides the exact posterior first and second central moments.  The Monte-Carlo based ensemble Kalman filter (EnKF) \cite{sr:evensen}, \cite{sr:evensen03} \cite{Evensen2000} has been shown to provide computationally robust posterior samples for high dimensional nonlinear applications, which has thus spurred the theoretical analysis of both the EnKf (i.e. in  discrete time) \cite{LeGland2009} \cite{TongKelly2016} \cite{Majda2016} and ensemble Kalman-Bucy filter (EnKBF) (i.e. in continuous time time) \cite{Bergemann2012}, \cite{deWiljesStannat2018}, \cite{deWiljes2020}, \cite{Lange2019} \cite{Schillings2016} \cite{Bishop_P_2018} \cite{DelMoral2017}.  These references are by no means exhaustive and the literature on this topic continues to expand.  Despite its desirable computational properties particularly in the small ensemble setting, the EnKF and EnKBF do not produce samples from the true posterior, even as the number of particles $N \rightarrow \infty$ due to the linear update approximation.  Traditional sequential Monte Carlo approaches (also known as particle filters) are asymptotically consistent but have had limited practical application due to know well-known weight degeneracy issues even in problems of moderate dimension (e.g. $d < 100$).  Over the last decade, there has been increasing attention given to asymptotically consistent controlled filtering methods in the like of \cite{sr:crisan10}, \cite{sr:reich10} and in particular, the feedback particle filter (FPF) \cite{Yang2013}, \cite{Laugesen2015}, \cite{Yangmfc2011}.  These methods involve constructing a McKean-Vlasov process to guide the evolution of the particles in state space, thereby combining the robustness advantages of the EnKBF and the consistency properties of traditional SMC.  The process is constructed to ensure that the conditional distribution of the process at time $t$ given $\mathcal{Z}_t$ coincides exactly with the solution of the Kushner-Stratonovich equation.  The general form of this process is \cite{Pathiraja2021}
	\begin{equation}
		\begin{aligned}
			\label{eq:generalmckean}
			d\tilde{X}_t &= \mathcal{M}(\tilde{X}_t)dt + d\tilde{V}_t + K(\tilde{X}_t, \tilde{\rho}_t) \left( dZ_t -\frac{1}{2}(h(\tilde{X}_t) + \hat{\tilde{h}}_t)dt \right) \\
			&+ \frac{1}{2}(\nabla K(\tilde{X}_t, \tilde{\rho}_t) )^T K(\tilde{X}_t, \tilde{\rho}_t)dt + \mathcal{J}(\tilde{X}_t, \tilde{\rho}_t) dt,
		\end{aligned}
	\end{equation}
	where $\mathcal{J}$ and $K$ satisfy
	\begin{align}
		\nabla \cdot (\tilde{\rho}_t \mathcal{J}(x, \tilde{\rho}_t)) &= 0, \\
		\label{eq:Kfirstorderpde}
		\nabla \cdot (\tilde{\rho}_t K(x, \tilde{\rho}_t)) &= - (h - \hat{\tilde{h}}_t) \tilde{\rho}_t,
	\end{align}
	respectively, $\tilde{\rho}_t$ indicates the conditional density of the process given $\mathcal{Z}_t$ and $\hat{\tilde{h}}_t = \int h(x) \tilde{\rho}_t(x) dx$.  The non-uniqueness of solutions of (\ref{eq:Kfirstorderpde}) gives rise to many formulations established in the literature, and can be further analysed to obtain new formulations with desirable analytic and numerical properties.  A further unifying analysis of these filters can be found in \cite{Pathiraja2021}.  The focus of this paper is the FPF, a special case of (\ref{eq:generalmckean}) which can be obtained for the choice that $K$ is of gradient type.  More specifically, the exact FPF equations are
	\begin{equation}
		\label{eq:exactFPF}
		\begin{aligned}
			d\tilde{X}_t^F &= \mathcal{M}(\tilde{X}_t^F)dt + d\tilde{V}_t^F + K(\tilde{X}_t^F, {\rho}_t) \left( dZ_t -\frac{1}{2}(h(\tilde{X}_t^F) + \hat{{h}}_t)dt \right) \\
			&+ \frac{1}{2}(\nabla K(\tilde{X}_t^F, {\rho}_t) )^T K(\tilde{X}_t^F, {\rho}_t)dt,
		\end{aligned}
	\end{equation}
	where $K = \nabla \phi$ and $\phi$ satisfies the weighted Poisson equation
	\begin{align}
		\label{eq:poissexact}
		-\Delta_{\rho_t} \phi(x):= -\frac{1}{\rho_t}\nabla \cdot (\rho_t \nabla \phi(x) = (h - \hat{h}_t),
	\end{align}
	for strictly positive densities $\rho_t$, where $\rho_t$ refers to the conditional density of the process $\tilde{X}_t^F$.  $K$ is typically referred to as the gain function, due to its similarity with the Kalman gain in the Kalman filter.  The well-posedness of (\ref{eq:poissexact}) under general conditions on $\mathcal{M}$ and $h$ in unbounded domains is a largely open research question, with some first steps in \cite{Pathiraja2021} and \cite{Laugesen2015}.  A computational bottleneck of this approach is in obtaining an approximate solution to (\ref{eq:poissexact}).  An attractive feature of the FPF is that adopting the so-called constant gain approximation in (\ref{eq:exactFPF}) leads to the enKBf.  Various methods have been proposed for higher order approximations; see for instance using Galerkin based PDE solution methods with POD adjustment for high dimensional applications see e.g. \cite{Berntorp2016} and also using dynammic programming \cite{Radhakrishnan2016} and also more recently based on the active area of using neural networks for approximating solutions of high dimensional PDEs \cite{Olmez2020}.  Here we focus on the recently proposed method of \cite{Taghvaei2019} which makes use of diffusion maps.  Diffusion maps \cite{Coifman2006} is a manifold learning technique which uses the eigenvectors of the diffusion operator on a data set to determine a typically lower dimensional underlying manifold from which the data has been sampled.  This is done by first constructing a family of anisotropic diffusion processes which are the limits of a graph Laplacian jump process.  This family is parameterised by $\alpha$ which reflects the influence of the density of the data points in determining the infinitessimal transitions of the diffusion.  The key point here is to extend the classical ideas of using the normalised graph Laplacian to cluster data for a graph to the Euclidean space setting by making use of the arrangement or density of the data points.  Alternative normalisations are possible and advantageous, see e.g. an approach based on Sinkhorn weights in \cite{Wormell2021}. For further details on diffusion maps, see in particular the seminal work of \cite{Coifman2006}.

	\subsection{Diffusion map based approximation of the gain}
	The aforementioned properties of diffusion maps make it an ideal candidate for developing a sample based approximation of $\phi$.  \cite{Taghvaei2019} uses the diffusion map approach in a novel way by starting from the semigroup formulation of the Poisson equation (\ref{eq:poissexact}).  Their algorithm is summarised in the remainder of this section, we refer to \cite{Taghvaei2019} for further details and derivations.  The semigroup formulation of (\ref{eq:poissexact}) is given by
	\begin{align}
		\label{eq:fixedpteqninfinite}
		\phi = P_\epsilon \phi + \int_0^\epsilon P_s(h - \hat{h})ds,
	\end{align}
	where $\{P_\epsilon\}_{\epsilon \geq 0}$ is the Markov semigroup associated to the weighted Laplacian $\Delta_{\rho_t}$, which must be approximated.  In other words, $P_s$ is the transition operator of a Langevin diffusion with invariant density given by $\rho_t$.  Adopting $\alpha = \frac{1}{2}$ in the diffusion map method corresponds to a Markov chain approximating the Langevin diffusion of interest \cite{Coifman2006}, with anisotropic transition kernel denoted by $k_\epsilon$ and
	\begin{align*}
		k_\epsilon(x,y) := \frac{\tilde{g}_\epsilon(x,y)}{(\int \tilde{g}_\epsilon(x,z) \rho_t(z)dz)^{1/2} (\int \tilde{g}_\epsilon(y,z) \rho_t(z)dz)^{1/2}},
	\end{align*}
	where
	\begin{align}
		\label{eq:gaussiankern}
		\tilde{g}_\epsilon(x,y) := (4 \pi \epsilon)^{-d/2} \exp(-\frac{1}{4 \epsilon} |x-y|^2),
	\end{align}
	is a Gaussian kernel with bandwidth matrix $2\epsilon I_d$ where $I_d$ denotes the identity matrix of size $d$.  This can be used to obtain the following approximation of the transition operator $P_\epsilon$, denoted by ${T}_\epsilon$,
	\begin{align}
		\label{eq:diffmapoperator}
		{T}_\epsilon f(x) := \int p_\epsilon(x,y)  f(y) dy,
	\end{align}
	where
	\begin{align*}
		p_\epsilon(x,y) := \frac{k_\epsilon(x,y)  \rho_t(y)}{k_\epsilon * \rho_t (x)},
	\end{align*}
	where we use the shorthand notation $k_\epsilon * \rho_t(x) := \int k_\epsilon(x,y) \rho_t(y) dy$, also in the remainder of the article.
	Additionally, we will use the following equivalent simplified form throughout the article
	\begin{align}
		\label{eq:peps}
		p_\epsilon(x,y) = \frac{g_\epsilon(x,y) \nu_\epsilon(y)}{g_\epsilon * \nu_\epsilon(x)},
	\end{align}
	where $g_\epsilon(x,y):= \exp(-\frac{1}{4 \epsilon}\norm{x-y}^2)$ is the unnormalised Gaussian kernel and
	\begin{align}
		\label{eq:defnnu}
		\nu_\epsilon(y) = \frac{\rho_t(y)}{(g_\epsilon * \rho_t(y))^{1/2}}.
	\end{align}
	An infinite dimensional, diffusion map based approximation of $\phi$, denoted by $\phi_\epsilon$, is then obtained as a solution of
	\begin{align}
		\label{eq:phidiffmap}
		\phi_\epsilon =  {T}_\epsilon \phi_\epsilon + \epsilon(h - \hat{h}_\epsilon); \quad \hat{h}_\epsilon := \int h(x) \pi_\epsilon (x) dx,
	\end{align}
	where $\pi_\epsilon(x):= \frac{\rho_t(x) k_\epsilon * \rho_t(x)}{\int \rho_t(x) k_\epsilon * \rho_t(x) dx} $ is the invariant density associated with the Markov operator ${T}_\epsilon$.  Now given $N$ $\rho_t$-distributed $\mathbb{R}^{d \times 1}$-valued random samples denoted by $\{ X^i\}_{i \in \{1, 2, \cdots, N\}}$, a finite dimensional approximation of ${T}_\epsilon$ is constructed as an $N \times N$ Markov matrix denoted by $T$ with $(i,j)$th entry given by
	\begin{subequations}
		\begin{align}
			\label{eq:tmp}
			T_{ij} & := \frac{k_{\epsilon}^N(X^i,X^j)}{\sum_l {k}_{ \epsilon}^N(X^i, X^l)}, \\
			\label{eq:Tij}
			& =  \frac{\tilde{q}_{\epsilon}(X^i,X^j)}{{\sum_l \tilde{q}_{ \epsilon}(X^i, X^l)}},
		\end{align}
	\end{subequations}
	where
	\begin{align*}
		\quad {k}_{ \epsilon}^N(X^i, X^j) & := \frac{\tilde{g}_\epsilon(X^i, X^j)}{(\sum_l \tilde{g}_\epsilon(X^j, X^l))^{1/2}(\sum_l \tilde{g}_\epsilon(X^i, X^l))^{1/2}}, \\
		\tilde{q}_{ \epsilon}(X^i, X^j) &:= \frac{g_\epsilon(X^i, X^j)}{(\sum_l g_\epsilon(X^j, X^l))^{1/2}}.
	\end{align*}
	In the remainder of the article, we will adopt the shorthand notation
	\begin{align}
		\label{eq:defnseps}
		s_\epsilon^i := \sum_l \tilde{q}_\epsilon (X^i, X^l).
	\end{align}
	Furthermore, $T$ is a reversible Markov matrix with strictly positive entries and a unique stationary distribution (see proposition 4.1 \cite{Taghvaei2019}).  $T$ is then used to construct a finite dimensional approximation of the fixed point equation (\ref{eq:phidiffmap}),
	\begin{align}
		\label{eq:fixedpteqn}
		\Phi = T \Phi + \epsilon ({\bf h} - \hat{h}^N \mathbf{1}),
	\end{align}
	where $\Phi \in \mathbb{R}^{N \times 1}$,  ${\bf h}:= [h(X^{1}), h(X^{2}), \cdots, h(X^{N})]^T$, $\mathbf{1} \in \mathbb{R}^{N \times 1}$ is a vector of 1's and $\hat{h}^N = \sum_{i=1}^N \pi_i h(X^i)$ where $\pi_i = \frac{\sum_l \tilde{q}_{ \epsilon}(X^i, X^l)}{\sum_j \sum_l \tilde{q}_{ \epsilon}(X^j, X^l)}$ is the stationary distribution of the Markov matrix $T$.  The solution of (\ref{eq:fixedpteqn}) is denoted by $\Phi^\infty$ and it holds that
	\begin{align}
		\label{eq:phidiff}
		\Phi^{\infty}_j - \Phi^{\infty}_i = \epsilon (h(X^j) - h(X^i)) +  \epsilon \sum_{n=1}^{\infty} [T^n  {\bf h}]_j -   [T^n  {\bf h}]_i .
	\end{align}
	Taking the derivative of an intermediate approximation of (\ref{eq:phidiffmap}) given by
	\begin{align}
		\label{eq:phiinterm}
		\phi_\epsilon^{(N)}(x) := \frac{1}{\sum_i \tilde{q}_\epsilon(x, X^i)} \sum_{j=1}^N \tilde{q}_\epsilon(x, X^j) \Phi_j + \epsilon(h(x) - \hat{h}^N),
	\end{align}
	(see Appendix C in \cite{Taghvaei2019}) leads to the following finite dimensional approximation of the gain $K = \nabla \phi$
	\begin{align*}
		K_\epsilon(X^i, \rho^N)  = \frac{1}{2\epsilon} \sum_j \sum_k T_{ij} T_{ik}(r_j - r_k) X^j,
	\end{align*}
	where
	\begin{align}
		\label{eq:rjdefn}
		r_j := \Phi_j^\infty + \epsilon h(X^j),
	\end{align}
	and $\rho^N$ denotes the empirical measure of the particles; which can be equivalently expressed as
	\begin{subequations}
		\begin{align}
			\label{eq:tmp2}
			K_\epsilon(X^i, \rho^N)  &=  \sum_{j=1}^N s_{ij} X^j, \\
			\label{eq:KfiniteN}
			&:= \frac{1}{2 \epsilon} \sum_j T_{ij} X^j r_j - \frac{1}{2 \epsilon} \left(\sum_j T_{ij} X^j \right) \left(\sum_{j} T_{ij} r_j\right),
		\end{align}
	\end{subequations}
	which has the benefit of making explicit the weighted cross-covariance type structure, thereby showing the resemblance to the Kalman gain in the enKBf.  An attractive feature of the diffusion map based approximation of the gain is that it leads to a bridging between the constant gain approximation (i.e. the enKBf) and the exact filter via the $\epsilon$ parameter.  More specifically, as $\epsilon \rightarrow \infty$, the approximation (\ref{eq:KfiniteN}) approaches the ensemble Kalman gain matrix.  For a precise explanation, see \cite{Taghvaei2019}.

	\subsection{Problem statement}
	Throughout the article, we consider the following $N$ particle process approximating (\ref{eq:exactFPF}) excluding the Stratonovich term, with the diffusion map based approximation of $K$ from \cite{Taghvaei2019},
	\begin{align}
		\label{eq:finiteNprocess}
		dX_t^{i} = \mathcal{M}(X_t^{i})dt + dV_t^i + K_\epsilon(X_t^{i}, \rho_t^N) \left( dZ_t - \frac{1}{2} (h(X_t^{i}) + \hat{h}_t^N) dt   \right),
	\end{align}
	where $\rho_t^N$ denotes the empirical measure of the $N$ particle process at time $t$ and $K_\epsilon(X_t^{i}, \rho_t^N) $ is a finite dimensional diffusion map based approximation of $K$ defined in (\ref{eq:KfiniteN}).  The Stratonovich term is neglected in (\ref{eq:finiteNprocess}) primarily to simplify the analysis; the influence of this term in practical applications will generally be quite small for large $N$ and with Lipschitz $\mathcal{M}$ and $h$ where the distributions are close to the Gaussian setting (since the gain $K$ will be close to constant in $x$ so that the derivative term is close to zero.)
	
	The main goal of this article is two-fold; 1) to determine the well-posedness of the finite $N$ particle system (\ref{eq:finiteNprocess}); and 2) to establish the existence of a mean-field limit of (\ref{eq:finiteNprocess}) as $N \rightarrow \infty$ for fixed $\epsilon > 0$.  In order to address the second point, we first propose that the limiting form of (\ref{eq:finiteNprocess}) can be described by the following McKean-Vlasov process
	\begin{align}
		\label{eq:meanfield}
		d\bar{X}_t = \mathcal{M}(\bar{X}_t)dt + d\bar{V}_t + \bar{K_\epsilon}(\bar{X}_t, \bar{\rho}_t) \left( dZ_t - \frac{1}{2} (h(\bar{X}_t) + \hat{\bar{h}}_t) dt   \right),
	\end{align}
	where $\bar{\rho}_t$ indicates the density function associated to the conditional law of the process given the observations up to time $t$ and $\hat{\bar{h}}_t = \int h(x) \rho_t(x) dx$.  The gain function $\bar{K}_\epsilon$ is the mean-field form of (\ref{eq:KfiniteN}) and is defined as
	\begin{align*}
		\bar{K}_\epsilon(x, \bar{\rho}_t) := \frac{1}{2\epsilon} \int \mathcal{R}_\epsilon(y) y  p_\epsilon(x,y) dy - \frac{1}{2\epsilon} \left(  \int y p_\epsilon(x,y) dy\right)   \left( \int \mathcal{R}_\epsilon(y)  p_\epsilon(x,y) dy  \right),
	\end{align*}	
	where $\mathcal{R}_\epsilon(y):= \phi_\epsilon(y) + \epsilon h(y)$ and $\phi_\epsilon(y)$ is the solution of the fixed point equation (\ref{eq:phidiffmap}) and the transition density $p_\epsilon(x,y)$ is as defined in (\ref{eq:peps}), but with $\rho_t$ replaced by $\bar{\rho}_t.$
	
	The well-posedness of the diffusion map approximation of $\phi$ as well as the convergence in both $\epsilon \rightarrow 0$ (bias) and $N \rightarrow 0$ (variance) was investigated in \cite{Taghvaei2019} for the case of iid samples from a density of the form $\rho(x) = \exp(-V(x))$ with $V(x) = \frac{1}{2}(x - m)^T \Sigma (x-m) + w(x)$ where $w(x)$ is smooth and uniformly bounded.  The convergence of $\phi_\epsilon^{(N)}$ as defined in (\ref{eq:phiinterm}) to $\phi_\epsilon$ as $N \rightarrow \infty$ was established for densities supported on compact domains.
	
	\subsection{Statement of the main results}
	\label{sec:main}
	
	\begin{theorem}
		\label{theo:wellposedfinite}
		\textbf{Well-posedness of the finite $N$ system.} Assume $\mathcal{M}: \mathbb{R}^d \rightarrow \mathbb{R}^d$ is globally Lipschitz with Lipschitz constant $L_M$ and $h: \mathbb{R}^d \rightarrow \mathbb{R}$ is uniformly bounded and globally Lipschitz with Lipschitz constant $L_h$.  Then the $N$-approximate interacting particle system (\ref{eq:finiteNprocess}) with finitely many iterates used to solve (\ref{eq:fixedpteqn}) posesses a unique strong solution for all $t > 0$.
	\end{theorem}	
	
	The proof of Theorem \ref{theo:wellposedfinite} is given in Section \ref{sec:wellposedfiniteN}.  As will be made clear there, the well-posedness result does not hold uniformly in $N$ and well-posedness of the mean field process (\ref{eq:meanfield}) is non-trivial due to the fact that the coefficients have non-Lipschitz dependence on the measure argument.   However, the properties of $\bar{K}_\epsilon$ with respect to the spatial argument established in Section \ref{sec:analmeanK} indicate a promising pathway to well-posedness for processes defined on unbounded domains, assuming that the conditional density of the process has sufficiently fast tail decay.  We capture this more precisely via log-concavity, see Assumption \ref{ass:logconc}.  The next main result relates to the convergence of the finite $N$ process to its mean field limit.  Propagation of chaos type results have been established for the enKBf in the linear case \cite{DelMoral2017} and also more general ensemble square root filters \cite{Lange2020}.  The covariance structure of the gain lends itself well to such analysis.  The weighted covariance structure of (\ref{eq:KfiniteN})leads to several complications; our main result is therefore a trajectorial propagation of chaos up to a stopping time.  That is,
	\begin{theorem}
		\label{theo:propchaos}
		\textbf{Trajectorial propagation of chaos.}
		Given $h$ and $\mathcal{M}$ globally Lipschitz with Lipschitz constants $L_h$ and $L_M$ respectively, $h$ uniformly bounded, and for (\ref{eq:fixedpteqn}) solved exactly and $\epsilon > 0$, the following mean square error convergence result holds up to a stopping time $\zeta_\delta^N > 0$ for large enough $N$,
		\begin{align}
			\label{eq:l2control}
			\mathbb{E} \left[ \sup_{t \in [0,\zeta_\delta^N]} \frac{1}{N}\sum_{i=1}^N \left| X_t^{i} - \bar{X}_t^i \right|^2  \right]  \leq {\frac{C}{N}},
		\end{align}
		where
		\begin{align*}
			C = f(\epsilon, d, L_M, L_h, \norm{h}_\infty),
		\end{align*}
		and $X_t^{i}$ is the $i$th particle solution of (\ref{eq:finiteNprocess}) with iid initial conditions $\{X_0^{i,N}\}_{i \in \{1, 2, \cdots N \}}$ and $\{\bar{X}_t^i\}_{i \in \{1, 2, \cdots N \}}$ are $N$ independent solutions of (\ref{eq:meanfield}) with the same Brownian motions $\{V_t^i\}$ and  initial conditions used to generate $\{{X}_t^{i}\}_{i \in \{1, 2, \cdots N \}}$.
	\end{theorem}	
	The precise form of the stopping time and conditions for large enough $N$ are given in Section \ref{sec:poc}, along with the proof of the above theorem in Section \ref{sec:proofthm2}.  The stopping time arguments used in the proof of this theorem allow us to consider the exact solution to (\ref{eq:fixedpteqn}), unlike the $L$-approximate solution considered in Theorem \ref{theo:wellposedfinite}.  The precise form of the constant $C$ can also be found in Section \ref{sec:proofthm2}, where it can be seen that it scales rather poorly with the dimension $d$, although this maybe used to specify a criterion for choosing $\epsilon$.  The complicated structure of the diffusion map approximation to the gain function means that standard propagation of chaos techniques for interacting particle systems e.g. \cite{Bolley2010} are not directly applicable here, leading to a non-trivial analysis.  The analysis also differs from that used to obtain standard stability results in discrete time particle filters.  There are several open questions which should be explored in future research, in addition to the well-posedness mentioned above.  In particular, extending Theorem \ref{theo:propchaos} to arbitrary time intervals $[0,T]$ for any $T > 0$ by way of apriori estimates on $\bar{X}$ and $X$ to ensure non-explosion in finite time (i.e. a finite time well-posedness result).  As discussed above, the well-posedness of the mean field process remains yet to be established and is an open problem in the literature \cite{Lange2020}.  Such localisation arguments usually allow for convergence in probability, although $L^2$ convergence as in Theorem \ref{theo:wellposedfinite} may be possible using similar arguments via Moore Osgood theorem as used in \cite{P_wongzakaiarxiv2020}.  Finally, it would be of course pertinent to consider the gain function approximation in the exact filter (\ref{eq:exactFPF}) and to consider the simultaneous limits $\epsilon \rightarrow 0$ and $N \rightarrow \infty$.  The precise relation between $\epsilon$ and $N$ would most likely need to be carefully controlled, as indicated by the bias variance discussion in \cite{Taghvaei2019}.  We note that the convergence of the fixed point equation on unbounded domains is still an open question in general\cite{Taghvaei2019} and would also need to be established.  On a more practical note, it may also be worth considering variable bandwidth kernels for which several advantages have been established, see e.g. \cite{Berry2016}.
	
	The remainder of the article is structured as follows.  Section \ref{sec:analmeanK} provides some analytic properties of the mean field gain $K_\epsilon(x,\mu)$ subject to certain conditions on $\mu$.  Section \ref{sec:wellposedfiniteN} examines well-posedness of the $N$ interacting particle equations for fixed $N$, i.e. the proof of the first main result.  Finally, Section \ref{sec:poc} contains the proof of second main result along with several important lemmas necessary for proving Theorem \ref{theo:propchaos}.

	\section{Analytic properties of $\bar{K}_\epsilon(x, \mu)$}
	\label{sec:analmeanK}
	This section presents some insights on properties of the mean field form of the gain, $\bar{K}_\epsilon$, in unbounded domains by making use of a log-concavity assumption.  They are mainly of independent interest and maybe utilised for obtaining propagation of chaos type results as in Section \ref{sec:poc} without the use of stopping times.  The exponential tail decay of the conditional density ensured by the aforementioned assumption plays an important role in uniformly bounding the gain in the spatial argument, which is a rather remarkable result given the covariance type structure of the gain and without boundedness assumptions on $h$.  More specifically, in this section we assume the following.  Note that we use the notation $\succcurlyeq$ to indicate a matrix inequality, i.e. $A \succcurlyeq B$ means $(A-B)$ is positive semidefinite.
	\begin{aspt}
		\label{ass:logconc}
		The conditional density of the mean field process (\ref{eq:meanfield}) at time $t$ is strongly log-concave, i.e. it takes the form $\exp(-\mathcal{V}_t(x))$ with
		\begin{align*}
			\nabla_x^2 \mathcal{V}_t \succcurlyeq  c_v I,
		\end{align*}
		for some $c_v > 0$.	
	\end{aspt}

	\begin{lemma}
		\label{lem:logconc}
		\textbf{Log-concavity of the transition kernel} Under Assumption \ref{ass:logconc}, it holds that the transition kernel $p_\epsilon(x,y)$at time $t$ is strongly log-concave in $y$ for all $x \in \mathbb{R}^d$.  That is, $p_\epsilon(x,y) \propto \exp(-\mathcal{G}_t(x,y))$ with
		\begin{align*}
			\nabla_y^2 \mathcal{G}_t \succcurlyeq c_g I,
		\end{align*}
		where
		\begin{align}
			\label{eq:cgdefn}
			c_g = \frac{1}{4 \epsilon} + c_v.
		\end{align}
		\begin{proof}
			The unnormalised form of $p_\epsilon (x,y)$ can be rewritten as
			\begin{equation*}
				\begin{aligned}
					g_\epsilon(x,y) \nu_\epsilon(y) &= \exp \left(-\frac{1}{4 \epsilon} |x-y|^2 - \mathcal{V}_t(y) - \frac{1}{2} \log(g_\epsilon * \bar{\rho}_t(y)) \right) \\
					& =: \exp(-\mathcal{G}_t (x,y)),
				\end{aligned}
			\end{equation*}
			and
			\begin{align*}
				\nabla_y^2 \mathcal{G}_t  = \frac{1}{2\epsilon} I + \nabla_y^2 \mathcal{V}_t(y) + \frac{1}{2} \left( \frac{1}{4 \epsilon^2} Cov_{\omega(y,z)dz}(z)  -\frac{1}{2\epsilon} I \right),
			\end{align*}
			where $\omega(y,z) \propto \exp(-F_t(y,z))$ and $F_t(y,z) = \frac{1}{4 \epsilon} |y-z|^2 + \mathcal{V}_t(z)$.  Differentiation under the integral sign is justified by the same arguments as in the proof of Theorem 4.2 in \cite{Brascamp1976}.  Due to the positive-definiteness of $Cov_{\omega(y,z)dz}(z)$, it holds that
			\begin{align*}
				\nabla_y^2 \mathcal{G}_t     & \succcurlyeq \frac{1}{2 \epsilon} I + c_vI - \frac{1}{4 \epsilon} I,
			\end{align*}
			which gives the result.		
		\end{proof}	
	\end{lemma}
	
	The above lemma guarantees several nice properties of the gain function $\bar{K}_\epsilon$, as shown in the following lemmas.  In particular, it ensures Lipschitz continuity of $K_\epsilon$ in the spatial argument with a Lipschitz constant depending on the spectral gap of $p_\epsilon(x,dy)$.
	
	\begin{lemma}
		\label{lem:boundKbar}
		\textbf{Uniform boundedness of $\bar{K}_\epsilon(x,\mu)$ and $\nabla_x \bar{K}_\epsilon (x, \mu)$} in $x$.  Let $\mu$ be any log-concave density as per Assumption \ref{ass:logconc} with $c_v > \frac{1}{4 \epsilon}$and $h:\mathbb{R}^d \rightarrow \mathbb{R}$ a $C^1$ and globally Lipschitz function.  Then the following uniform bounds hold
		\begin{align}
			\label{eq:Kepsbound}
			\norm{\bar{K}_\epsilon(x,\mu)} &\leq \frac{\norm{\nabla h}_\infty }{2 c_g} \left(\frac{4\epsilon c_g -1}{2 \epsilon c_g - 1} \right)  \quad \forall \enskip x \in \mathbb{R}^d \\
			\label{eq:gradKepsbound}
			\norm{\nabla_x \bar{K}_\epsilon(x,\mu) }_F &\leq \frac{d \norm{\nabla h}_\infty}{ \epsilon} \frac{1}{c_g^{3/2}}  \left(\frac{4\epsilon c_g -1}{2 \epsilon c_g - 1} \right),  \quad \forall \enskip x \in \mathbb{R}^d
		\end{align}
		where $c_g$ is the spectral gap of $p_\epsilon(x,y)$, as defined in (\ref{eq:cgdefn}) and $\norm{\cdot}_F$ indicates the Frobenius norm.
		
		\begin{proof}
			
			Consider $\phi_\epsilon$, the solution of (\ref{eq:phidiffmap}).  For small enough $\epsilon$, the Neumann series expansion of (\ref{eq:phidiffmap}) gives
			\begin{subequations}
				\begin{align}
					\label{eq:tmp3}
					\phi_\epsilon^\infty (y) & = \epsilon(h(y) - \hat{h}_\epsilon) + \epsilon \sum_{n =1}^\infty \int_{E^{n}}   (h(y_n) - \hat{h}_\epsilon) \prod_{k=1}^n p_\epsilon (y_{k-1}, y_k)  dy_1 dy_2 \cdots dy_{n} \\
					\label{eq:phineumann}
					& = \epsilon \sum_{n=0}^\infty T_\epsilon^n (h(y) - \hat{h}_\epsilon),
				\end{align}		
			\end{subequations}
			where $y_0 = y$ and $\int_{E^n} $ indicates $n$ iterated integrals. The Markov operator $T_\epsilon h(y) = \int h(z) p_\epsilon (y,z) dz$ is globally Lipschitz in $y$ since for any $v \in \mathbb{R}^d$
			\begin{subequations}
				\begin{align}
					\label{eq:tmp4}
					\langle \nabla_y T_\epsilon h(y), \nabla_y T_\epsilon h(y) \rangle  & = \frac{1}{2 \epsilon} \int h(z) \langle \nabla_y T_\epsilon h(y) , z - p_\epsilon * z(y) \rangle p_\epsilon(y,z) dz \\
					\label{eq:tmp5}
					& = \frac{1}{2 \epsilon} \int (h(z) - \hat{h}_p) \langle \nabla_y T_\epsilon h(y) , z - p_\epsilon * z(y)  \rangle p_\epsilon(y,z) dz \\
					\label{eq:gradTbound}
					& \leq \frac{1}{2 \epsilon c_g}  \left( \int \norm{\nabla h}^2 p_\epsilon (y,z) dz  \right)^{1/2}  \norm{\nabla_y T_\epsilon h(y)},
				\end{align}  	
			\end{subequations}
			where the last inequality arises from Lemma \ref{lem:logconc} and a direct application of Theorem 1.1 in \cite{Carlen2013} (see also Lemma \ref{lem:covPI}).  This implies that
			\begin{align*}
				\norm{\nabla_y T_\epsilon h(y)}   & \leq \frac{\norm{\nabla h}_\infty }{2 \epsilon c_g},
			\end{align*}
			which combined with (\ref{eq:phineumann}) implies that $\phi_\epsilon^\infty$ is globally Lipschitz in the spatial argument, when $\mu$ satisfies Assumption \ref{ass:logconc} and $h$ is globally Lipschitz, since
			\begin{subequations}
				\begin{align}
					\label{eq:tmp6}
					\norm{\phi_\epsilon^\infty(x) - \phi_\epsilon^\infty(y)} & \leq   \epsilon \sum_{n=0}^\infty \norm{T^n_\epsilon h(x) - T^n_\epsilon h(y)} \\
					\label{eq:tmp7}
					& \leq \epsilon \sum_{n=0}^\infty  \left( \frac{1}{2\epsilon c_g} \right)^n  \norm{\nabla h}_\infty \norm{x - y} \\
					\label{eq:philipcons}
					& =  \left( \frac{\epsilon \norm{\nabla h}_\infty}{1 - \frac{1}{2 \epsilon c_g}} \right) \norm{x-y},
				\end{align}
			\end{subequations}
			which holds so long as $2\epsilon c_g >1$, which is guaranteed by the assumption on $c_v$.  By the same calculations as in (\ref{eq:gradTbound}), we have
			\begin{align*}
				\norm{\bar{K}_\epsilon(x, \mu)}^2 &= \frac{1}{2 \epsilon}  \int (\mathcal{R}_\epsilon(y) - p_\epsilon * \mathcal{R}_\epsilon(x)) \langle \bar{K}_\epsilon(x, \mu), y- p_\epsilon *y(x)  \rangle p_\epsilon (x,y) dy\\
				& \leq \frac{1}{2 \epsilon c_g} \left( \int \norm{\nabla \mathcal{R}_\epsilon}^2 p_\epsilon (x,y) dy  \right)^{1/2}   \norm{\bar{K}_\epsilon (x,\mu)}.
			\end{align*}
			Using (\ref{eq:philipcons}), we have that
			\begin{align}
				\label{eq:lipconstg}
				\norm{\nabla \mathcal{R}_\epsilon}_\infty \leq \epsilon \norm{\nabla h}_\infty \left(\frac{4\epsilon c_g -1}{2 \epsilon c_g - 1} \right),
			\end{align}
			which then gives the first result (\ref{eq:Kepsbound}).  In order to prove the second result, start with
			\begin{align*}
				\nabla_x \bar{K}_\epsilon(x,\mu) = -\frac{1}{4 \epsilon^2} \underbrace{\int  (\mathcal{R}_\epsilon(y) - p_\epsilon * \mathcal{R}_\epsilon(x))  (y - p_\epsilon * y(x)) (y - p_\epsilon * y(x))^T p_\epsilon (x,y) dy.}_{=: A}
			\end{align*}
			The Frobenius norm of $\nabla_x \bar{K}_\epsilon$ satisfies the inequality
			\begin{align}
				\label{eq:normgradK}
				\norm{\nabla_x \bar{K}_\epsilon(x,\mu) }_F \leq \frac{1}{4 \epsilon^2} \left(\sum_j \sum_i A_{ij}^2 \right)^{1/2}.
			\end{align}
			A uniform bound on the entries $A_{ij}$ can be obtained using similar calculations as above involving the Cauchy-Schwarz inequality and Lemmas \ref{lem:logconc} and \ref{lem:covPI}.  More specifically,
			\begin{align*}
				& \norm{A_{ij}} = \norm{\int (\mathcal{R}_\epsilon(y) - p_\epsilon * \mathcal{R}_\epsilon(x)) (y - p_\epsilon *y(x))_i (y - p_\epsilon *y(x))_j p_\epsilon (x,y) dy} \\
				& \leq \left(\int  (\mathcal{R}_\epsilon(y) - p_\epsilon * \mathcal{R}_\epsilon(x))^2 p_\epsilon (x,y) dy  \right)^{1/2}   \\
				& \times \left(\int   (y - p_\epsilon *y(x))_i^2 (y - p_\epsilon *y(x))_j^2 p_\epsilon (x,y) dy  \right)^{1/2} \\
				& \leq \frac{\norm{\nabla \mathcal{R}_\epsilon}_\infty}{\sqrt{c_g}} \left(\int   (y - p_\epsilon *y(x))_i^2 (y - p_\epsilon *y(x))_j^2 p_\epsilon (x,y) dy  \right)^{1/2} \\
				%			\label{eq:Aijbound}
				& \leq \frac{\norm{\nabla \mathcal{R}_\epsilon}_\infty}{\sqrt{c_g}} \left(\int   (y - p_\epsilon *y(x))_i^4 p_\epsilon (x,y) dy  \right)^{1/4} \left(\int   (y - p_\epsilon *y(x))_j^4 p_\epsilon (x,y) dy \right)^{1/4}.
			\end{align*}
			
			To bound the remaining terms, first let $f(x,y) := (y - p_\epsilon * y(x))^2$ and $f_i$ refer to the $i$th entry of the vector valued function $f$, so that
			\begin{align*}
				\int   (y - p_\epsilon *y(x))_i^4 p_\epsilon (x,y) dy &= \int   (f_i(x,y) - p_\epsilon * f_i(x) + p_\epsilon *f_i(x)   )^2 p_\epsilon (x,y) dy \\
				& \leq 2 \int (f_i(x,y) - p_\epsilon * f_i(x))^2 p_\epsilon (x,y) dy + 2 \left(p_\epsilon * f_i (x)  \right)^2.
			\end{align*}
			Applying Lemmas \ref{lem:logconc} and \ref{lem:covPI} to the two terms in the last equality yields
			\begin{align*}
				\int   (y - p_\epsilon *y(x))_i^4 p_\epsilon (x,y) dy  \leq \frac{8}{c_g^2} + \frac{2}{c_g^2} = \frac{10}{c_g^2}.
			\end{align*}
			Inserting the above estimate into the upper bound for $\norm{A_{ij}}$ and combining with (\ref{eq:normgradK})  and (\ref{eq:lipconstg}) yields the second result (\ref{eq:gradKepsbound}).
		\end{proof}	
	\end{lemma}

	\section{Well-posedness for the $N$-approximate process}
	\label{sec:wellposedfiniteN}
	Let {\small ${\bf X}_t = [X_t^1; X_t^2; ...; X_t^N]\\ \in \mathbb{R}^{dN \times 1}$} be the vector of particles, and also define {\small ${\bf V}_t = [V_t^1; V_t^2; ...; V_t^N] \in \mathbb{R}^{(dN+1) \times 1}$}.  We can rewrite the $N$-FPF in the following form
	\begin{align}
		\label{eq:mckeanvlasovmultivariateform}
		d{\bf X}_t = B({\bf X}_t) dt + I_{dN} d{\bf V}_t + \mathcal{K}_\epsilon({\bf X}_t) dZ_t,
	\end{align}
	where $I_{dN}$ indicates the identity matrix of size $dN$, $B: \mathbb{R}^{dN} \rightarrow \mathbb{R}^{dN}$ and $\mathcal{K}_\epsilon: \mathbb{R}^{dN} \rightarrow \mathbb{R}^{dN}$ are measurable functions given by
	\begin{align*}
		B({\bf X}_t) :=  \begin{bmatrix}
			b(X^1_t, L_N({\bf X}_t)) \\
			b(X_t^2, L_N({\bf X}_t)) \\
			\vdots \\
			b(X_t^N, L_N({\bf X}_t))
		\end{bmatrix};	\quad
		\mathcal{K}_\epsilon({\bf X}_t) := \begin{bmatrix}
			K_\epsilon(X^1_t, L_N({\bf X}_t)) \\
			K_\epsilon(X_t^2, L_N({\bf X}_t)) \\
			\vdots \\
			K_\epsilon(X_t^N, L_N({\bf X}_t)),
		\end{bmatrix}
	\end{align*}
	and
	\begin{align*}
		b(X_t^i, L_N({\bf X}_t)) &:= \mathcal{M}(X_t^i) - \frac{1}{2} K_\epsilon(X_t^i, L_N({\bf X}_t)) (h(X_t^i) + \hat{h}_t^N),
	\end{align*}
	where $\rho_t^N = L_N({\bf X}_t)$ and $L_N: (\mathbb{R}^d)^N \rightarrow \mathcal{P}(\mathbb{R}^d)$ is the empirical measure map given by \medskip
	\begin{align*}
		L_N(x_1, x_2, \cdots, x_N) = \frac{1}{N} \sum_{k=1}^N \delta_{x_k}.
	\end{align*}
	It then follows from standard Ito theory that if $B,\mathcal{K}_\epsilon$ are $C^2$ (and therefore locally Lipschitz) and of linear growth, that the SDE (\ref{eq:mckeanvlasovmultivariateform}) has a unique strong solution. The proof of Theorem \ref{theo:wellposedfinite} follows immediately from the following lemma and the fact that $h$ is assumed to be uniformly bounded and $\mathcal{M}$ globally Lipschitz.  As will be clear in the proof of Lemma \ref{lem:Klienargrowth}, this well-posedness only holds for finite $N$.  Furthermore, we restrict our attention to cases where the fixed point equation is only solved for $L$ iterates, whose solution we denote by $\Phi^{(L)}$.  This is due to the difficulty in obtaining a control on the mixing rate of the Markov chain uniformly in the particles without apriori conditions, see e.g. the stopping time criterion in Section \ref{sec:poc}.  In the following, we drop the time subscript as the results hold for all $t \geq 0$.
	
	\begin{lemma}
		\label{lem:Klienargrowth}
		Assume that $h:\mathbb{R}^d \rightarrow \mathbb{R}$ is globally Lipschitz and uniformly bounded and the fixed point equation (\ref{eq:fixedpteqn}) solved for $L$ iterates.  Then the function $\mathcal{K}_\epsilon(\bf X)$ satisfies
		\begin{align}
			\label{eq:Klineargrowth}
			\norm{\mathcal{K}_\epsilon(\bf X)} < C_k(1 + \norm{{\bf X}}),
		\end{align}
		where
		\begin{align*}
			C_k = N^6  d^{3/2} \sqrt{\epsilon} (1 \vee \epsilon) (L_h + L \norm{h}_\infty).
		\end{align*}
		\begin{proof}
			Recall the definition of $K_\epsilon$ from (\ref{eq:KfiniteN}). We first show that $s_{ij}$ is uniformly bounded for all $i,j$.  Specifically,
			\begin{align*}
				\norm{s_{ij}} &= \left| \frac{1}{2 \epsilon} T_{ij} \left( \sum_{k=1}^N T_{ik} (r_j  - r_k) \right) \right| \\
				\nonumber
				& \leq \frac{1}{2 \epsilon}  \left( \sum_{k=1}^N T_{ij} T_{ik} |r_j - r_i| + T_{ij} T_{ik} |r_i  - r_k| \right) \\
				& \leq \frac{1}{2 \epsilon}  \left(  T_{ij} |r_j - r_i| +  \sum_{k=1}^N T_{ik} |r_i  - r_k| \right).
			\end{align*}
			For the terms in the summand, notice that
			\begin{align}
				\label{eq:Tikrirkupperbound}
				T_{ik}|r_i  - r_k| & \leq T_{ik}| \Phi_i^{(L)} - \Phi_k^{(L)}| + \epsilon T_{ik}L_h|X^i -  X^k |.
			\end{align}
			It follows from Jensen's inequality that $T_{ik} \leq N g_\epsilon(X^i, X^k)$  and therefore
			\begin{align}
				\label{eq:Tikxixkbound}
				T_{ik}|X^i -  X^k | \leq  N 2\sqrt{d \epsilon} \exp\left(-\frac{1}{2} \right).
			\end{align}
			Additionally, we have that
			\begin{align*}
				\norm{\Phi_k^{(L)} - \Phi_i^{(L)}} &\leq \epsilon L_h\norm{X^k - X^i} +  \epsilon \norm{\sum_{n=1}^{L} [T^n  {\bf h}]_k -   [T^n  {\bf h}]_i} \\
				& \leq \epsilon L_h\norm{X^k - X^i} +  \epsilon \norm{h}_\infty L N \sum_j \norm{T_{kj} - T_{lm}}.
			\end{align*}
			
			The summand can be analysed starting from (\ref{eq:Tij}) and recalling the definition of $s_\epsilon^i$ in (\ref{eq:defnseps}) as 
			\begin{align*}
				\norm{T_{kj} - T_{ij}} &\leq \frac{1}{s_\epsilon^k} \norm{\tilde{q}_\epsilon (X^k, X^j)  - \tilde{q}_\epsilon (X^i, X^j)} + \tilde{q}_\epsilon (X^i, X^j) \frac{1}{s_\epsilon^k s_\epsilon^i} \norm{s_\epsilon^i - s_\epsilon^k} \\
				&  < \frac{1}{s_\epsilon^k} \norm{\tilde{q}_\epsilon (X^k, X^j)  - \tilde{q}_\epsilon (X^i, X^j)} +  \frac{1}{s_\epsilon^k s_\epsilon^i} \norm{s_\epsilon^i - s_\epsilon^k} \\
				& < \sqrt{N} \norm{\tilde{q}_\epsilon (X^k, X^j)  - \tilde{q}_\epsilon (X^i, X^j)} +  N \norm{s_\epsilon^i - s_\epsilon^k} \\
				& < \sqrt{N} \norm{g_\epsilon (X^k, X^j)  - g_\epsilon (X^i, X^j)} +  N \sum_l \norm{g_\epsilon (X^k, X^l)  - g_\epsilon (X^i, X^l)}  \\
				& < 2N^2 K_g \norm{X^k - X^i},
			\end{align*}
			where for the last inequality, we have used (\ref{eq:lipg}).
			Combining all gives
			\begin{align*}
				\norm{\Phi_k^{(L)} - \Phi_i^{(L)}} &< D_\phi \norm{X^k - X^i},
			\end{align*}
			where $D_\phi: = \epsilon \left( L_h  +  4 LN^3 \norm{h}_\infty d \epsilon \exp(-1) \right)$.  Substituting the above into (\ref{eq:Tikrirkupperbound}) gives
			\begin{align*}
				T_{ik}|r_i  - r_k| & \lesssim N^4 (d \epsilon)^{3/2} (1 \vee \epsilon) (L_h + L \norm{h}_\infty),
			\end{align*}
			so that for all $i,j$,
			\begin{align*}
				\norm{s_{ij}} \lesssim \underbrace{  N^5 d^{3/2} \sqrt{\epsilon} (1 \vee \epsilon) (L_h + L \norm{h}_\infty).}_{=: C_s}
			\end{align*}
			The result then immediately follows from the fact that
			\begin{align*}
				\norm{\mathcal{K}_\epsilon({\bf X})} < N C_s \norm{\bf X}.
			\end{align*}
		\end{proof}
	\end{lemma}

	\section{Propagation of chaos}
	\label{sec:poc}
	Here we utilise the conditional propagation of chaos technique to establish the mean field limit of (\ref{eq:finiteNprocess}).  The standard coupling technique involves considering $N$ independent copies of (\ref{eq:meanfield}), $\{\bar{X}_t^i\}_{i \in \{1, 2, \cdots N \}}$ with the same Brownian motions $\{V^i_t\}$ and iid initial conditions $\{X_0^i\}$, i.e.
	\begin{align*}
		d\bar{X}_t^i = \mathcal{M}(\bar{X}_t^i)dt + dV_t^i + \bar{K_\epsilon}(\bar{X}_t^i, \bar{\rho}_t) \left( dZ_t - \frac{1}{2} (h(\bar{X}_t^i) + \hat{\bar{h}}_t) dt   \right); \quad \bar{X}^i_0 = X_0^i,
	\end{align*}
	where we use $\bar{\rho}_t$ with a slight abuse of notation.  We also use $\rho_t^N$ to denote the empirical measure of the particles.  The conditional case simply involves conditioning on the common noise, in this case, the observation process $Z$.

	As mentioned in Section \ref{sec:main}, we provide a result up to a stopping time only due to the difficulties in establishing convergence of $\bar{K}_\epsilon(\bar{X}_t^i, \bar{\rho}_t)$ to $K_\epsilon(X_t^i, \rho_t^N)$.  The following stopping time arguments effectively restrict the results to compact domains.  Furthermore, this ensures that there exists a unique solution $\phi_\epsilon$ as established formally in \cite{Taghvaei2019} for a measure $\mu$ with bounded support (and satisfying further conditions).  Well-posedness of the finite $N$ process follows from the same arguments as in the proof of Theorem \ref{theo:wellposedfinite}; the bounded $h$ assumption is easily removed when considering a time interval up to the stopping time.  The overall goal is to first establish a uniform in $N$ control on the term $\norm{K_\epsilon(X_t^i, \rho_t^N) - \bar{K}_\epsilon(\bar{X}_t^i, \bar{\rho}_t)}$ which is facilitated by several intermediate results of a similar nature involving the solution of the fixed point equation (Lemma \ref{lem:philip} and \ref{lem:rrbarlip}) and on the transition matrix (Lemma \ref{lem:lipT}).  The final crucial ingredient is a strong law of large numbers type result, given in Lemma \ref{ProofSLLN}.
	
	More specifically, consider the stopping times
	\begin{align*}
		\tau_\delta^N &:= \inf \{ t \geq 0:  \min_{1 \leq i \leq N} \frac{1}{N} \sum_k \tilde{q}_{\epsilon}(X^i_t, X^k_t) < \delta  \} \\
		\bar{\tau}_\delta^N &:=  \inf \{ t \geq 0:  \min_{1 \leq i \leq N}  \frac{1}{N} \sum_k \tilde{q}_{\epsilon}(\bar{X}^i_t, \bar{X}^k_t) < \delta \} \\
		\bar{\tau}_\delta & :=  \inf \{ t \geq 0:  \min_{1 \leq i \leq N} \int g_\epsilon(\bar{X}^i_t,y) \rho_t(y) dy < \delta  \} \\
	\end{align*}
	where $0 < \delta < 1$.   Furthermore, define
	\begin{align*}
		\zeta_\delta^N := T \wedge \tau_\delta^N \wedge \bar{\tau}_\delta^N  \wedge \bar{\tau}_\delta
	\end{align*}
	
	\subsection{Intermediate results}
	Before proving Theorem \ref{theo:propchaos}, we establish some important intermediate results.  All results here apply for $t < \zeta_\delta^N$. Clearly we have
	\begin{align}
		\label{eq:qtilecond}
		s_\epsilon^i := \sum_{k=1}^N \tilde{q}_{\epsilon}(X^i_t, X^k_t) \geq \delta N  \quad \forall \enskip i.
	\end{align}
	Additionally, it holds that
	\begin{align*}
		\tilde{q}_\epsilon (X_t^i, X_t^k) < g_\epsilon(X_t^i, X_t^k).
	\end{align*}
	This leads to the following upper and lower bounds (alternative characterisations are possible)
	\begin{align}
		\label{eq:Tbounds}
		\frac{1}{N^{3/2}} g_\epsilon(X_t^i,X_t^k) &< \frac{1}{N}\tilde{q}_{ \epsilon}(X_t^i, X_t^k) < T_{ik} < \frac{1}{\delta N} \tilde{q}_{ \epsilon}(X_t^i, X_t^k) < \frac{1}{\delta N} g_{ \epsilon}(X_t^i, X_t^k).
	\end{align}
	Furthermore, we have using (\ref{eq:qtilecond}) and (\ref{eq:Tbounds}) that
	\begin{align}
		\label{eq:boundq}
		\tilde{q}_\epsilon(X^j, X^m) &\leq \frac{g_\epsilon(X^j, X^m)}{\sqrt{\delta N}} \leq \frac{1}{\sqrt{\delta N}} \\
		\label{eq:boundT32}
		T_{jm} &\leq \frac{1}{(\delta N)^{3/2}}  \quad \forall \enskip j,m.
	\end{align}
	
	The aforementioned results also hold for these quantities evaluated on $\{\bar{X}_t^i\}_{i=1:N}$, e.g. where $\bar{T}_{ik}$ is given by
	\begin{align*}
		\bar{T}_{ik} := \frac{\tilde{q}_\epsilon(\bar{X}_t^i, \bar{X}_t^k)}{\sum_l \tilde{q}_\epsilon(\bar{X}_t^i, \bar{X}_t^l)}.
	\end{align*}
	We first introduce the following condition on $N$ which aids in obtaining estimates on the fixed point equation (\ref{eq:phidiff}).  Many of the estimates described below are therefore valid for ``large enough'' $N$ depending on the choice of $\delta$, specifically,
	\begin{aspt}
		\label{ass:delt}
		For a given $0 < \delta <1$,
		\begin{align*}
			N > \frac{1}{4 \delta^3}.
		\end{align*}
	\end{aspt}

	The first intermediate result is the following Lipschitz type property, its proof can be found in Appendix \ref{app:lipT}.  Throughout this section, we use the notation $a \lesssim b$ to mean $a \leq Cb$ where $C$ is a constant independent of $\epsilon, d, N, \delta$.
	
	\begin{lemma}
		\label{lem:lipT}
		
		It holds for $t < \zeta_\delta^N$,  for all $i,k$ that 
		\begin{align}
			\label{eq:Tcontcond}
			\norm{T_{ik} - \bar{T}_{ik}}  &  \lesssim  d \epsilon \left(  \frac{1}{\delta^3 N^{3/2}} \norm{X^i - \bar{X}^i}  +   \frac{1}{\delta^2 N^{3/2}}  \norm{X^k - \bar{X}^k}  +   \frac{1}{\delta^4 N^2} \sum_l  \norm{X^l - \bar{X}^l}  \right)
		\end{align}
		and
		\begin{align}
			\label{eq:sumlipTsquared}
			\left( \sum_k \norm{T_{ik} - \bar{T}_{ik}} \right)^2 \lesssim d^2 \epsilon^2 \left(  \frac{1}{\delta^6 N} \norm{X^i - \bar{X}^i}^2  +   \frac{1}{\delta^8 N} \sum_l  \norm{X^l - \bar{X}^l}^2    \right)
		\end{align}
		as well as
		\begin{align}
			\nonumber
			& \left( \sum_k \norm{T_{ik}(X^k-X^i)  - \bar{T}_{ik}(\bar{X}^k-\bar{X}^i) } \right)^2 \\
			\label{eq:sumlipTXsquared}
			& \lesssim   ((d \epsilon)^3 \vee K_f^2)\left( \frac{1}{\delta^6 N}\norm{X^i - \bar{X}^i}^2 + \frac{1}{\delta^4 N}  \sum_l \norm{X^l - \bar{X}^l}^2    \right),
		\end{align}
		where $K_f$ is defined in (\ref{eq:Kfdefn}).	
	\end{lemma}	
	
	The following lemma shows that the difference of any two components of the fixed point equation (\ref{eq:fixedpteqn}) satisfy a kind of linear growth type property.  The proof can be found in Appendix \ref{app:lipphi} and makes use of a modification to the standard arguments for establishing contraction in the total variation metric for transition matrices with strictly positive entries.  The modified argument along with Condition (\ref{ass:delt}) allows us to indirectly satisfy a Doeblin minorization condition.
	
	\begin{lemma}
		\label{lem:philip}
		It holds for $t < \zeta_\delta^N$, with a fixed $0 < \delta < 1$, $h$ globally Lipschitz and uniformly bounded, that for all $k,l$ and $N$ satisfying condition (\ref{ass:delt}) that
		\begin{align}
			\label{eq:lipphi}
			\norm{r_k - r_l} \lesssim C_\phi \norm{X_t^k - X_t^l},
		\end{align}
		and
		\begin{align*}
			\norm{\bar{r}_k - \bar{r}_l} < C_\phi \norm{\bar{X}_t^k - \bar{X}_t^l},
		\end{align*}
		where
		\begin{align}
			\label{eq:cphi}
			C_\phi(N) &:=  \epsilon  \left( L_h +  \frac{{d \epsilon \norm{h}_\infty }}{2 \delta^{3/2}\sqrt{N} - 1}  \right),
		\end{align}
		and $\bar{r}_k$ is given by (\ref{eq:rjdefn}) evaluated on $\{ \bar{X}_t^i\}_{i=1:N}$.  Additionally, it holds that
		\begin{align}
			\label{eq:boundrdiff}
			\norm{r_k - r_l} \leq C_\gamma
		\end{align}
		where
		\begin{align*}
			C_\gamma := 2 \epsilon \norm{h}_\infty \left( \frac{5\delta^{3/2}\sqrt{N} - 2}{2 \delta^{3/2}\sqrt{N} - 1} \right),
		\end{align*}
		
		and likewise for $\bar{r}_k - \bar{r}_l$.	
	\end{lemma}	
	
	The contraction arguments used in the proof of the previous lemma also allow us to establish the following `Lipschitz' type property on $r$, analogous to Lemma \ref{lem:lipT}.  The proof can be found in Appendix \ref{app:rrbarlip}.
	\begin{lemma}
		\label{lem:rrbarlip}
		Assume the same conditions as Lemma \ref{lem:philip}.  Then it holds that for all $i,k$ and $N$ satisfying condition (\ref{ass:delt}) that
		\begin{align}
			\label{eq:rdiffLip}
			\norm{r_k - r_i - \bar{r}_i + \bar{r}_k}  & \lesssim   C_{r1} \left( \norm{X^i - \bar{X}^i} +\norm{X^k - \bar{X}^k} \right)  + C_{r2} \frac{1}{N} \sum_l  \norm{X^l - \bar{X}^l}
		\end{align}
		where
		\begin{align*}
			C_{r1} &:= \left( \epsilon L_h + d \epsilon^2 \norm{h}_\infty \left(\frac{ 1}{2 \delta^{3/2} \sqrt{N} - 1} \right) \frac{1}{\delta^{3/2} }  \right) \\
			C_{r2} & := \frac{\epsilon^2 d (L_h \vee \norm{h}_\infty)}{\delta^6}  \left(\frac{ 4 \delta^{3/2} \sqrt{N}-1}{2 \delta^{3/2} \sqrt{N} - 1} \right),
		\end{align*}
		from which it follows immediately that
		\begin{align}
			\label{eq:sumrdiffsquared}
			\sum_k \norm{r_k - r_i - \bar{r}_i + \bar{r}_k}^2  \lesssim C_{r1}^2 N \norm{X^i - \bar{X}^i}^2  +   (C_{r1}^2 + C_{r2}^2)  \sum_m \norm{X^m - \bar{X}^m}^2.
		\end{align}	
	\end{lemma}

	The following lemma establishes a couple of useful properties of the empirical $K_\epsilon$, its proof is given in Appendix \ref{app:Klip}.

	\begin{lemma}
		\label{lem:Klip}
		Assume the same conditions as Lemma \ref{lem:philip}.  Then for $N$ sufficiently large, i.e. satisfying condition (\ref{ass:delt}), we have the following bounds,
		\begin{align}
			\label{eq:Kuniform}
			\norm{K_\epsilon(X_t^i, \rho_t^N)} & \lesssim K_b,
		\end{align}
		where
		\begin{align*}
			K_b =  C_\phi  d,
		\end{align*}
		and
		\begin{align*}
			\frac{1}{N} \sum_{i=1}^N 	\norm{K_\epsilon(\bar{X}^i_t, \bar{\rho}_t^N) - K_\epsilon(X_t^i, \rho_t^N)}^2 & \lesssim C_K \frac{1}{N}  \sum_{i=1}^N \norm{X^i - \bar{X}^i}^2,
		\end{align*}
		where
		\begin{align*}
			C_K =   \frac{C_\gamma^2  ((d \epsilon)^3 \vee K_f^2)}{\delta^8 \epsilon^2} + \frac{d (C_{r1}^2 + C_{r2}^2) }{\epsilon}.
		\end{align*}
	\end{lemma}

	A consequence of the above Lemma is that $\nabla_{x} K_\epsilon(x, \mu_t^N)$ has a uniformly bounded derivative when $t < \zeta_\delta^N$ where the bound is independent of $N$ since $K_\epsilon(x, \mu_t^N)$ is continuously differentiable in $x$.  This will be useful when analysing the stratonovich correction term not considered here.

	The following lemma will be needed to obtain a LLN type bound and relies on the tail decay of the kernel $g_\epsilon(x,y)$ to control the functions $\phi$ and $X$ in $\bar{K}_\epsilon$.  We make use of the following assumptions
	\begin{aspt}
		\label{ass:existsoln}
		There exists a solution $\phi_\epsilon$ to the fixed point equation (\ref{eq:phidiffmap}).
	\end{aspt}	
	
	\begin{aspt}
		\label{ass:nubound}
		For $t < \zeta_\delta^N$, there exists a constant $\delta_v > 0$ such that
		\begin{align*}
			\int \norm{y} \nu_\epsilon(y) dy < \delta_v.
		\end{align*}
		where $\nu_\epsilon$ is defined in (\ref{eq:defnnu}).
	\end{aspt}	
	
	\begin{lemma}
		\label{ProofSLLN}
		Assume the same conditions as Lemma \ref{lem:philip}.  Then for $N$ sufficiently large, i.e. satisfying condition (\ref{ass:delt}) it holds that
		\begin{align}
			\label{eq:sllntest}
			\mathbb{E} \left[ \frac{1}{N} \sum_i \norm{ \bar{K}_\epsilon(\bar{X}_t^i, \bar{\rho}_t) - K_\epsilon(\bar{X}_t^i, \bar{\rho}_t^N)}^2 \right] \leq \frac{\mathcal{C}}{N},
		\end{align}
		where
		\begin{align*}
			\mathcal{C} = \frac{( (\norm{\phi_\epsilon}_\infty + \epsilon \norm{h}_\infty) \delta_v + C_\gamma d (\epsilon \vee \sqrt{\epsilon}) )^2}{\epsilon^2 \delta^4}.
		\end{align*}
	\end{lemma}
	The proof can be found in Appendix 	\ref{app:ProofSLLN}.

	\subsection{Proof of Theorem \ref{theo:propchaos}}
	\label{sec:proofthm2}
	We are now ready to prove Theorem \ref{theo:propchaos} related to convergence of the finite $N$ process to the mean field process up to a stopping time.  For any $i$, define
	\begin{align*}
		\mathcal{X}_t^i := \bar{X}_t^i - X_t^{i}.	
	\end{align*}	
	We have
	\begin{align*}
		&d\mathcal{X}_t^i =  [\mathcal{M}(\bar{X}_t^i) - \mathcal{M}(X_t^i) ] dt + \left[\bar{K}_\epsilon(\bar{X}_t^i,  \bar{\rho}_t) - K_\epsilon(X_t^i, \rho_t^N) \right] h(X_t^*) dt  \\
		& - \frac{1}{2} \left[ \bar{K}_\epsilon(\bar{X}_t^i,  \bar{\rho}_t) - K_\epsilon(X_t^i, \rho_t^N) \right] \left[ h(\bar{X}_t^i) + \widehat{\bar{h}}_t  \right] dt - \frac{1}{2} K_\epsilon(X_t^i, \rho_t^N) \left[ h(\bar{X}_t^i) - h(X_t^i)  \right]dt \\
		&- \frac{1}{2} K_\epsilon(X_t^i, \rho_t^N) \left[ \widehat{\bar{h}}_t - \hat{h}_t^N  \right]dt + [\bar{K}_\epsilon(\bar{X}_t^i, \bar{\rho}_t) - K_\epsilon(X_t^i, \bar{\rho}_t)] dW_t.
	\end{align*}
	Using (\ref{eq:Kuniform}) from Lemma \ref{lem:Klip}, we have that
	\begin{align*}
		&\sup_{s \in [0,t \wedge \zeta_\delta^N]}  \frac{1}{N} \sum_i \norm{\mathcal{X}_s^i}^2  \leq L_M^2   \int_0^{t \wedge \zeta_\delta^N} \sup_{u \in [0,s \wedge \zeta_\delta^N]}  \frac{1}{N} \sum_i \norm{\mathcal{X}_u^i}^2 ds    \\
		& + \norm{h}_\infty^2 \int_0^{t \wedge \zeta_\delta^N} \frac{1}{N} \sum_i \norm{\bar{K}_\epsilon(\bar{X}_s^i, \bar{\rho}_s) - K_\epsilon(X_s^i, \rho_s^N) }^2 ds \\
		&+ \frac{1}{4} L_h^2 K_b^2 \int_0^{t \wedge \zeta_\delta^N}  \sup_{u \in [0,s \wedge \zeta_\delta^N]} \frac{1}{N} \sum_i \norm{\mathcal{X}_u^i}^2 ds + \frac{1}{4}  K_b^2  \int_0^{t \wedge \zeta_\delta^N}   \norm{\widehat{\bar{h}}_s - \hat{h}_s^N }^2 ds \\
		& +  \sup_{s \in [0,t \wedge \zeta_\delta^N]} \frac{1}{N} \sum_i  \norm{\int_0^s  \bar{K}_\epsilon(\bar{X}_u^i, \bar{\rho}_u) - K_\epsilon(X_u^i, \rho_u^N) dW_u }^2 \\
		& := I_3 + I_4 + I_5 + I_6 + I_7,
	\end{align*}
	where $X_t^*$ denotes the references process.  Starting with the term $I_7$, it follows from Doob's maximal inequality and Ito isometry that
	\begin{align*}
		&\mathbb{E} \left[  I_7\right]  \lesssim \mathbb{E} \left[ \int_0^{t \wedge \zeta_\delta^N}  \frac{1}{N} \sum_i  \norm{ \bar{K}_\epsilon(\bar{X}_u^i, \bar{\rho}_u) - K_\epsilon(X_u^i, \rho_u^N)}^2 du \right].
	\end{align*}
	Furthermore, Lemma \ref{lem:Klip} and Lemma \ref{ProofSLLN} yield \allowdisplaybreaks
	\begin{align*}
		& \mathbb{E} \left[ \int_0^{t \wedge \zeta_\delta^N}  \frac{1}{N} \sum_i  \norm{ \bar{K}_\epsilon(\bar{X}_u^i, \bar{\rho}_u) - K_\epsilon(X_u^i, \rho_u^N)}^2 du \right] \\
		& \leq \int_0^{t \wedge \zeta_\delta^N}  \mathbb{E} \left[ \frac{1}{N} \sum_i \norm{ \bar{K}_\epsilon(\bar{X}_u^i, \bar{\rho}_u) - K_\epsilon(\bar{X}_u^i, \bar{\rho}_u^N)}^2 \right] du   \\
		&  +  \int_0^{t \wedge \zeta_\delta^N}  \mathbb{E} \left[\sup_{u \in [0,s \wedge \zeta_\delta^N]}  \frac{1}{N} \sum_i \norm{ {K}_\epsilon(\bar{X}_u^i, \bar{\rho}_u^N) - K_\epsilon(X_u^i, \rho_u^N)}^2 \right]  ds  \\
		& \leq \frac{\mathcal{C}}{N} (t \wedge \zeta_\delta^N)  + C_K \int_0^{t \wedge \zeta_\delta^N}  \mathbb{E} \left[ \sup_{u \in [0,s \wedge \zeta_\delta^N]} \frac{1}{N} \sum_i \norm{\mathcal{X}_u^i}^2 \right] ds,
	\end{align*}
	which can also be used to bound $I_4$.  To analyse $I_6$, notice that
	\begin{align*}
		\norm{\widehat{\bar{h}}_s - \hat{h}_s^N}^2 & \leq 2\norm{\widehat{\bar{h}}_s - \frac{1}{N} \sum_i h(\bar{X}_s^i)}^2 + 2\norm{\frac{1}{N} \sum_i h(\bar{X}_s^i) - h(X_s^i)}^2 \\
		& \leq 2\norm{\widehat{\bar{h}}_s - \frac{1}{N} \sum_j h(\bar{X}_s^i)}^2 + 2\frac{L_h^2}{N} \sum_i \norm{\mathcal{X}_s^i}^2,
	\end{align*}
	so that from the strong law of large numbers, we have
	\begin{align*}
		\mathbb{E} \left[   \int_0^{t \wedge \zeta_\delta^N}   \norm{\widehat{\bar{h}}_s - \hat{h}_s^N }^2 ds     \right] & \lesssim \frac{1}{N} \int_0^{t \wedge \zeta_\delta^N} Var(h(\bar{X}_s)) ds  \\
		& + L_h^2  \int_0^{t \wedge \zeta_\delta^N}  \mathbb{E} \left[ \sup_{u \in [0,s \wedge \zeta_\delta^N]}  \frac{1}{N} \sum_i \norm{\mathcal{X}_u^i}^2 \right]  ds  \\
		& \lesssim \frac{\norm{h}^2_\infty (t \wedge \zeta_\delta^N)}{N}  \\
		& + L_h^2  \int_0^{t \wedge \zeta_\delta^N}  \mathbb{E} \left[ \sup_{u \in [0,s \wedge \zeta_\delta^N]}  \frac{1}{N} \sum_i \norm{\mathcal{X}_u^i}^2 \right]  ds.
	\end{align*}
	Combining the above estimates along with Gronwall lemma gives
	\begin{align*}
		\mathbb{E} \left[  \sup_{s \in [0,t \wedge \zeta_\delta^N]}  \frac{1}{N} \sum_i \norm{\mathcal{X}_s^i}^2 \right]  \lesssim \frac{D_9\exp(D_8 (t \wedge \zeta_\delta^N))}{N},
	\end{align*}
	where
	\begin{align*}
		D_8 &:= (L_M^2 + L_h^2 K_b^2  + \norm{h}_\infty^2 C_K)  \\
		D_9 &:=  (K_b^2 + \mathcal{C} ) \norm{h}_\infty^2  (t \wedge \zeta_\delta^N).
	\end{align*}

	%For acknowledgements section, please don't number the section, please begin it with \section*{Acknowledgements}
	\section*{Acknowledgments} We are grateful for constructive comments from anonymous referees that have helped improve this article.

	\appendix
	
	\section{Useful analytic estimates}
	
	The following lipschitz property of $g_\epsilon(x,y)$ will be used extensively.
	\begin{align}
		\label{eq:lipg}
		\norm{g_\epsilon(x_1,y_1) - g_\epsilon(x_2, y_2)} \leq K_g \norm{x_1 - y_1 - x_2 + y_2},
	\end{align}
	with $K_g:= 2d\epsilon \exp(-1)$.  Additionally, we have the following result.
	
	\begin{lemma}
		\label{eq:flip}
		Consider $f: \mathbb{R}^d \rightarrow \mathbb{R}^d$ where $f(z) = z\exp(-\frac{1}{4 \epsilon} \norm{z}^2)$ for $\epsilon > 0$.  It holds that $f$ is globally lipschitz, i.e.
		\begin{align*}
			\norm{f(z) - f(y)} \leq K_f \norm{z - y},
		\end{align*}
		where
		\begin{align}
			\label{eq:Kfdefn}
			K_f = \sqrt{\exp(-2)d^2 + d(1 - \exp(-2))}.
		\end{align}
		\begin{proof}
			We start with $f_j(z) = z_j \exp(-\frac{1}{4 \epsilon} \norm{z}^2)$ which is a scalar function.  Then
			\begin{align*}
				\norm{\nabla_j f_j} &  = \frac{1}{2\epsilon}\exp(-\frac{1}{4 \epsilon} \norm{z}^2) \norm{2\epsilon - z_j^2} \\
				& \leq \frac{1}{2\epsilon}\exp(-\frac{1}{4 \epsilon} z_j^2) \norm{2\epsilon - z_j^2} \\
				& \leq 1,
			\end{align*}
			and also
			\begin{align*}
				\norm{\nabla_i f_j} & = \frac{1}{2 \epsilon} \norm{z_j z_i} \exp(-\frac{1}{4 \epsilon} \norm{z}^2), \quad i \neq j   \\
				& \leq \frac{1}{2 \epsilon} \norm{z_i} \exp(-\frac{1}{4 \epsilon}z_i^2) \norm{z_j}\exp(-\frac{1}{4 \epsilon}z_j^2) \\
				& \leq \exp(-1).
			\end{align*}
			So we have shown that the individual entries of the vector function $f$ are lipschitz, meaning
			\begin{align*}
				&\norm{f(z) - f(y)}\\& = \left( \sum_j (f_j(y) - f_j(z))^2 \right)^{1/2} \leq \left( \norm{y-z}^2 \sum_j  \left[ (\nabla_j f_j)^2 + \sum_{i \neq j} (\nabla_i f_j)^2 \right] \right)^{1/2} \\
				& \leq  \sqrt{\exp(-2)d^2 + d(1 - \exp(-2))} \norm{y-z}.
			\end{align*}		
		\end{proof}		
	\end{lemma}	
	
	\section{Proof of lemma \ref{lem:lipT}}
	\label{app:lipT}
	For notational ease, we drop the time subscript in the following which is held fixed.  We have that
	\begin{align*}
		T_{ik} - \bar{T}_{ik} = \frac{\tilde{q}_{ \epsilon}(X^i,X^k)}{\sum_l \tilde{q}_{ \epsilon}(X^i, X^l)} - \frac{\tilde{q}_{ \epsilon}(\bar{X}^i,\bar{X}^k)}{\sum_l \tilde{q}_{ \epsilon}(\bar{X}^i, \bar{X}^l)},
	\end{align*}
	and let us adopt the shorthand notation $\bar{s}_\epsilon^i = \sum_l \tilde{q}_{ \epsilon}(\bar{X}^i, \bar{X}^l).$  Then
	\begin{align*}
		\norm{T_{ik} - \bar{T}_{ik}} \leq \frac{\norm{\bar{s}^i_\epsilon - s^i_\epsilon}}{s^i_\epsilon \bar{s}^i_\epsilon} \tilde{q}_{ \epsilon}(X^i,X^k) + \frac{1}{\bar{s}_\epsilon^i} \norm{\tilde{q}_{ \epsilon}(X^i,X^k) - \tilde{q}_{ \epsilon}(\bar{X}^i,\bar{X}^k)  }.
	\end{align*}
	Similarly, we have using the shorthand notation $u_\epsilon^k := (\sum_l g_{ \epsilon}(X^k, X^l))^{1/2}$ that \allowdisplaybreaks
	\begin{align*}
		\norm{\tilde{q}_{ \epsilon}(X^i,X^k) - \tilde{q}_{ \epsilon}(\bar{X}^i,\bar{X}^k)} &= \norm{\frac{g_{ \epsilon}(X^i,X^k)}{u_\epsilon^k} - \frac{g_{ \epsilon}(\bar{X}^i,\bar{X}^k)}{\bar{u}_\epsilon^k} } \\
		& \leq  \frac{\norm{\bar{u}^k_\epsilon - u^k_\epsilon}}{u^k_\epsilon \bar{u}^k_\epsilon} g_{ \epsilon}(X^i,X^k) + \frac{1}{\bar{u}_\epsilon^k} \norm{g_{ \epsilon}(X^i,X^k) - g_{ \epsilon}(\bar{X}^i,\bar{X}^k)  } \\
		& \leq \frac{\norm{\bar{u}^k_\epsilon - u^k_\epsilon}}{u^k_\epsilon \bar{u}^k_\epsilon} + \frac{1}{\bar{u}_\epsilon^k} K_g \left(\norm{X^i - \bar{X}^i}  + \norm{X^k - \bar{X}^k} \right)  \\
		& \leq \frac{\norm{\bar{u}^k_\epsilon - u^k_\epsilon}}{\delta N} + \frac{1}{\sqrt{\delta N}} K_g \left(\norm{X^i - \bar{X}^i}  + \norm{X^k - \bar{X}^k} \right).
	\end{align*}
	In the above we have used  (\ref{eq:Tbounds}) and (\ref{eq:lipg}).  Now notice that
	\begin{align*}
		u_\epsilon^k = \sqrt{\sum_l (g_{2 \epsilon}(X^k, X^l))^2} = \norm{G_{2 \epsilon}({\bf z})},
	\end{align*}
	where $G_{2 \epsilon}({\bf z}) := \left[g_{2 \epsilon}(z_1), g_{2 \epsilon}(z_2), \cdots, g_{2 \epsilon}(z_N)   \right]^T$ and $z_l:= X^k - X^l$.   By the reverse triangle inequality and  (\ref{eq:lipg}), we have
	\begin{align*}
		\norm{u_\epsilon^k - \bar{u}_\epsilon^k} & = \norm{\norm{G_{2 \epsilon}({\bf z})} - \norm{G_{2 \epsilon}({\bf \bar{z}})}}  \\
		&\leq \norm{G_{2 \epsilon}({\bf z}) - G_{2 \epsilon}({\bf \bar{z}})} \\
		& \leq 2 K_{g}\left(  \sum_l \norm{X^k - X^l - (\bar{X}^k -\bar{X}^l)}^2     \right)^{1/2}  \\
		&\hspace{-2cm} \leq 2\sqrt{2} K_g \left(\sqrt{N} \norm{X^k - \bar{X}^k} + \left( \sum_l  \norm{X^l - \bar{X}^l}^2 \right)^{1/2} \right).
	\end{align*}
	Combining gives
	\begin{align*}
		\norm{\tilde{q}_{ \epsilon}(X^i,X^k) - \tilde{q}_{ \epsilon}(\bar{X}^i,\bar{X}^k)} & \lesssim \frac{K_g}{ \sqrt{\delta N}} \norm{X^i - \bar{X}^i}  + K_g \left( \frac{1}{\delta \sqrt{N}} + \frac{1}{\sqrt{\delta N}} \right)  \norm{X^k - \bar{X}^k}  \\
		& + \frac{K_g}{\delta N} \left( \sum_l  \norm{X^l - \bar{X}^l}^2 \right)^{1/2},
	\end{align*}
	from which it follows that
	\begin{equation*}
		\begin{aligned}
			\norm{\bar{s}^i_\epsilon - s^i_\epsilon} \lesssim \frac{K_g \sqrt{N}}{\sqrt{\delta}}\norm{X^i - \bar{X}^i} & + K_g \left( \frac{1}{\delta \sqrt{N}} +  \frac{1}{\sqrt{\delta N}}    \right) \sum_l \norm{X^l - \bar{X}^l} \\
			\label{eq:sepsdiff}
			& + \frac{K_g}{\delta} \left( \sum_l  \norm{X^l - \bar{X}^l}^2 \right)^{1/2},
		\end{aligned}
	\end{equation*}
	and then using (\ref{eq:boundq}) and (\ref{eq:qtilecond}), we have that
	\begin{align*}
		\norm{T_{ik} - \bar{T}_{ik}}  \lesssim  d \epsilon \left(  \frac{1}{\delta^3 N^{3/2}} \norm{X^i - \bar{X}^i}  +   \frac{1}{\delta^2 N^{3/2}}  \norm{X^k - \bar{X}^k}  +   \frac{1}{\delta^4 N^2} \sum_l  \norm{X^l - \bar{X}^l}  \right),
	\end{align*}
	which gives the first result (\ref{eq:Tcontcond}). The result  (\ref{eq:sumlipTsquared}) follows in a straightforward manner using Jensen's inequality.
	The calculations to obtain (\ref{eq:sumlipTXsquared}) are entirely analogous (making use of Lemma \ref{eq:flip}  instead of (\ref{eq:lipg})) and are therefore omitted.
	\qed

	\section{Proof of lemma \ref{lem:philip}}
	\label{app:lipphi}
	
	Recall (\ref{eq:phidiff}).  For $n=0$, we have trivially that
	\begin{align}
		\label{eq:boundT0lip}
		\norm{[T^0  {\bf h}]_k -   [T^0  {\bf h}]_l }\leq L_h \norm{X^k - X^l},
	\end{align}
	and also
	\begin{align}
		\label{eq:boundT0h}
		\norm{[T^0  {\bf h}]_k -   [T^0  {\bf h}]_l }\leq 2 \norm{h}_\infty,
	\end{align}
	and for $n=1$ it holds that
	\begin{align*}
		\norm{[T  {\bf h}]_k -   [T  {\bf h}]_l}  \leq \norm{h}_\infty \sum_m \norm{T_{km} - T_{lm}}.
	\end{align*}
	Now consider $n=2$, where it holds that
	\begin{align*}
		\norm{[T^2  {\bf h}]_k -   [T^2  {\bf h}]_l} \leq \norm{h}_\infty \sum_m \norm{[T^2]_{km} - [T^2]_{lm}},
	\end{align*}
	and
	\begin{align*}
		\sum_m \norm{[T^2]_{km} - [T^2]_{lm}} &= \sum_m \norm{ \sum_j (T_{kj} - T_{lj}) T_{jm}} \\
		&= 2 \norm{(v_k - v_l)T}_{TV},
	\end{align*}
	where $v_k := [T_{k1}, T_{k2}, \cdots, T_{kN}]$.  Since $v_{k}$ and $v_{l}$ are probability vectors and $T$ is a Markov matrix with strictly positive entries, it holds that $T$ is a strict contraction in the total variation metric.  By a slight modification of standard arguments and making us of (\ref{eq:boundT32}), we can define a matrix $Q$ with entries
	\begin{align*}
		Q_{jm} = \frac{\delta\sqrt{\delta N}}{1-\delta\sqrt{\delta N}} \left( \frac{1}{(\delta N)^{3/2}} - T_{jm}  \right),
	\end{align*}
	which is a stochastic matrix since all the entries are non-negative and $\sum_m Q_{jm} = 1$.
	Therefore, we can write $T_{jm} = \frac{1}{(\delta N)^{3/2}} - \frac{1-\delta\sqrt{\delta N}}{\delta\sqrt{\delta N}} Q_{jm}$
	\begin{align*}
		&\sum_m \norm{ \sum_j (v_k(j) - v_l(j)) T_{jm}} \\ &= \sum_m \norm{ \sum_j (v_k(j) - v_l(j)) \left(  \frac{1}{(\delta N)^{3/2}} - \frac{1-\delta\sqrt{\delta N}}{\delta\sqrt{\delta N}} Q_{jm} \right) } \\
		& \leq \frac{1}{(\delta N)^{3/2}}\sum_m \norm{ \sum_j (v_k(j) - v_l(j))} \\
		& + \sum_m \norm{ \sum_j \left(  (v_l(j) - v_k(j)) \frac{1-\delta\sqrt{\delta N}}{\delta\sqrt{\delta N}} Q_{jm}  \right)}  \\
		& = \frac{1-\delta\sqrt{\delta N}}{\delta\sqrt{\delta N}} \sum_j  \norm{ v_l(j) - v_k(j)} \sum_m  Q_{jm} \\
		& = \frac{1-\delta\sqrt{\delta N}}{\delta\sqrt{\delta N}} 2 \norm{v_k - v_l}_{TV}.
	\end{align*}
	That is, we have that
	\begin{align*}
		\sum_m \norm{[T^2]_{km} - [T^2]_{lm}} & \leq \left(\frac{1}{\delta^{3/2} \sqrt{N}} - 1 \right) \sum_j \norm{T_{kj} - T_{lj}}.
	\end{align*}
	Inserting the above along with (\ref{eq:boundT0h}) into (\ref{eq:phidiff}) and induction gives
	\begin{align*}
		\norm{\Phi^{\infty}_k - \Phi^{\infty}_l } &=  \epsilon  \norm{\sum_{n=0}^{\infty} [T^n  {\bf h}]_k -   [T^n  {\bf h}]_l }  \\
		& \leq 2\epsilon \norm{h}_\infty  + \epsilon \norm{h}_\infty \left(\sum_j \norm{T_{kj} - T_{lj}} \right) \sum_{n=1}^\infty \left(\frac{1}{\delta^{3/2} \sqrt{N}} - 1 \right)^{n-1}  \\
		& = 2\epsilon \norm{h}_\infty + \epsilon \norm{h}_\infty \left(\sum_j \norm{T_{kj} - T_{lj}} \right) \frac{1}{2 - \frac{1}{\delta^{3/2} \sqrt{N}}},
	\end{align*}
	where the geometric series converges so long as $\norm{\frac{1}{\delta^{3/2} \sqrt{N}} - 1 } < 1$, which is satisfied under condition (\ref{ass:delt}).  The result (\ref{eq:boundrdiff}) then follows immediately.  To obtain the result (\ref{eq:lipphi}), use (\ref{eq:boundq}) to obtain
	\begin{align*}
		&\norm{T_{km} - T_{lm}} \leq \frac{1}{s^k_\epsilon} \norm{\tilde{q}_\epsilon (X^k, X^m) - \tilde{q}_\epsilon (X^l, X^m)} + \tilde{q}_\epsilon(X^l, X^m) \frac{1}{s_\epsilon^k s_\epsilon^l} \norm{s_\epsilon^l - s_\epsilon^k} \\
		& \leq  \frac{1}{\delta N} \norm{\tilde{q}_\epsilon (X^k, X^m) - \tilde{q}_\epsilon (X^l, X^m)} + \frac{1}{(\delta N)^{5/2}} \norm{s_\epsilon^l - s_\epsilon^k} \\
		& \leq \frac{1}{\delta N} \norm{\tilde{q}_\epsilon (X^k, X^m) - \tilde{q}_\epsilon (X^l, X^m)} + \frac{1}{(\delta N)^{5/2}} \sum_j \norm{\tilde{q}_\epsilon (X^k, X^j) - \tilde{q}_\epsilon (X^l, X^j)},
	\end{align*}
	and it holds that
	\begin{align*}
		\norm{\tilde{q}_\epsilon (X^k, X^m) - \tilde{q}_\epsilon (X^l, X^m)} &\leq \frac{1}{\sqrt{\sum_j g_\epsilon(X^m, X^j)}} \norm{g_\epsilon(X^k, X^m) - g_\epsilon(X^l, X^m)} \\
		& \leq \frac{1}{\sqrt{\delta N}} {K_g} \norm{X^k - X^l},
	\end{align*}
	using (\ref{eq:lipg}).  Combining gives
	\begin{align}
		\label{eq:sumTdiff}
		\sum_m \norm{T_{km} - T_{lm}} &\leq \frac{{2K_g}}{\delta^{3/2} \sqrt{N}}  \norm{X^k - X^l}.
	\end{align}
	Again, inserting the above along with (\ref{eq:boundT0lip}) into (\ref{eq:phidiff}) gives
	\begin{align*}
		\norm{\Phi^{\infty}_k - \Phi^{\infty}_l } & \leq \epsilon L_h \norm{X^k - X^l} + \epsilon \norm{h}_\infty \frac{{2K_g}}{\delta^{3/2} \sqrt{N}}  \norm{X^k - X^l} \sum_{n=1}^\infty  \left(\frac{1}{\delta^{3/2} \sqrt{N}} - 1 \right)^{n-1}  \\
		& \leq \epsilon L_h \norm{X^k - X^l} + \epsilon \norm{h}_\infty \frac{{4 \exp(-1) d \epsilon}}{\delta^{3/2} \sqrt{N}}  \norm{X^k - X^l} \frac{1}{2 - \frac{1}{\delta^{3/2} \sqrt{N}}}.
	\end{align*}
	\qed

	\section{Proof of lemma \ref{lem:rrbarlip}}
	\label{app:rrbarlip}
	Here we will make use of the following simplifying notation
	\begin{align*}
		\mathcal{T}_i &:= \sum_l \norm{T_{il} - \bar{T}_{il}},  \\
		\mathcal{G}_k &:= \sum_l T_{kl} \sum_m \norm{T_{lm} - \bar{T}_{lm}}.
	\end{align*}
	With a slight abuse of notaiton, we also use $\mathcal{G}_k$ to denote $\sum_l \bar{T}_{kl} \sum_m \norm{T_{lm} - \bar{T}_{lm}}$.  Now it holds that
	\begin{align*}
		&\norm{r_k - r_i - \bar{r}_i + \bar{r}_k}  \leq \epsilon \left( 2L_h \norm{X^k - \bar{X}^k} + 2L_h \norm{X^i - \bar{X}^i}   \right)  \\
		& + \epsilon \sum_{n=1}^\infty \norm{ \underbrace{[T^n ({\bf h} - \hat{h} {\bf 1})]_k - [T^n ({\bf h} - \hat{h} {\bf 1})]_i - [\bar{T}^n ({\bf \bar{h}} - \hat{\bar{h}} {\bf 1})]_k + [\bar{T}^n ({\bf \bar{h}} - \hat{\bar{h}} {\bf 1})]_i}}.
	\end{align*}
	We denote the summand in the above by $ \mathcal{S}_r^n$.  Notice that the summand can be written as
	\begin{align*}
		\mathcal{S}_r^n & = \sum_m ([T^n]_{km} - \pi_m - [\bar{T}^n]_{km} + \bar{\pi}_m)h(X^m)  \\
		& - \sum_m ([T^n]_{im} - \pi_m - [\bar{T}^n]_{im} + \bar{\pi}_m)h(X^m) \\
		&+ \sum_m ( [\bar{T}^n]_{km} - \bar{\pi}_m) (h(X^m) - h(\bar{X}^m)) -\sum_m ( [\bar{T}^n]_{im} - \bar{\pi}_m) (h(X^m) - h(\bar{X}^m)) \\
		& =: S_k^n -S_i^n + I_k^n - I_i^n.
	\end{align*}
	Starting with $S_k$,
	\begin{align*}
		S_k^n - S_i^n & = \sum_m h(X^m) \sum_l  ([T^{n-1}]_{kl} - [\bar{T}^{n-1}]_{kl})T_{lm} + ([\bar{T}^{n-1}]_{kl} - \bar{\pi}_l)(T_{lm} - \bar{T}_{lm} ) \\
		& - ([T^{n-1}]_{il} - [\bar{T}^{n-1}]_{il})T_{lm} - ([\bar{T}^{n-1}]_{il} - \bar{\pi}_l)(T_{lm} - \bar{T}_{lm}),  \\
	\end{align*}
	So that
	\begin{align*}
		\norm{S_k^n - S_i^n} & \leq \norm{h}_\infty \left(\sum_m \norm{\sum_l  ([T^{n-1}]_{kl} - [\bar{T}^{n-1}]_{kl})T_{lm}} \right. \\
		& \left. + \sum_m \norm{\sum_l ([T^{n-1}]_{il} - [\bar{T}^{n-1}]_{il})T_{lm}  } \right. \\
		& \left. + \sum_m  \norm{\sum_l ([\bar{T}^{n-1}]_{kl} - \bar{\pi}_l)(T_{lm} - \bar{T}_{lm} ) - ([\bar{T}^{n-1}]_{il} - \bar{\pi}_l)(T_{lm} - \bar{T}_{lm})} \right)  \\
		& =: \norm{h}_\infty(U_k^n + U_i^n+ V_{ik}^n).
	\end{align*}
	In order to analyse the second term in the above, consider first the case $n=1$.  Then we have
	\begin{align*}
		V_{ik}^1 & \leq \mathcal{T}_k + \mathcal{T}_{i}.
	\end{align*}
	Then for $n=2$
	\begin{align*}
		V_{ik}^2 &\leq \sum_m  \norm{\sum_l \bar{T}_{kl} (T_{lm} - \bar{T}_{lm} ) - \bar{T}_{il}(T_{lm} - \bar{T}_{lm})}  \\
		& \leq \mathcal{G}_k + \mathcal{G}_i.
	\end{align*}
	Then for $n=3$,
	\begin{align*}
		V_{ik}^3 &\leq  \sum_m \sum_l  \norm{([\bar{T}^2]_{kl} - \bar{\pi}_l)(T_{lm} - \bar{T}_{lm} )} + \norm{([\bar{T}^2]_{il} - \bar{\pi}_l)(T_{lm} - \bar{T}_{lm})}  \\
		& \leq \sum_j \left( \norm{\bar{T}_{kj} - \bar{\pi}_j } + \norm{\bar{T}_{ij} - \bar{\pi}_j }  \right) \mathcal{G}_j,
	\end{align*}
	and $n=4$,
	\begin{align*}
		V_{ik}^4 &\leq  \sum_m \sum_l  \norm{([\bar{T}^3]_{kl} - \bar{\pi}_l)(T_{lm} - \bar{T}_{lm} )} + \norm{([\bar{T}^3]_{il} - \bar{\pi}_l)(T_{lm} - \bar{T}_{lm})}  \\
		&  \leq \sum_j \left( \norm{[\bar{T}^2]_{kj} - \bar{\pi}_j } + \norm{[\bar{T}^2]_{ij} - \bar{\pi}_j }  \right) \mathcal{G}_j.
	\end{align*}
	Now note that for $n \geq 1$, the following contraction holds
	\begin{subequations}
		\begin{align}
			\label{eq:tmp8}
			\sum_j \norm{[\bar{T}^n]_{kj} - \bar{\pi}_j} & \leq  \left( \frac{1}{\delta^{3/2} \sqrt{N}} - 1  \right)^{n-1} \sum_j \norm{T_{kj} - \pi_j}  \\
			\label{eq:boundTtostat}
			& \leq 2 \left( \frac{1}{\delta^{3/2} \sqrt{N}} - 1  \right)^{n-1},
		\end{align}
	\end{subequations}
	by similar reasoning as in the proof of Lemma \ref{lem:philip}.
	So that by induction, for $n \geq 3$,
	\begin{align*}
		V_{ik}^n &  \leq \sum_j \left( \norm{[\bar{T}^{n-2}]_{kj} - \bar{\pi}_j } + \norm{[\bar{T}^{n-2}]_{ij} - \bar{\pi}_j }  \right) \mathcal{G}_j  \\
		& \leq 4 \left( \frac{1}{\delta^{3/2} \sqrt{N}} - 1  \right)^{n-3} {\mathcal{G}},
	\end{align*}
	where $\mathcal{G}:= \max_j \mathcal{G}_j$.  Finally,
	\begin{align*}
		\sum_{n=1}^\infty V_{ik}^n & \lesssim \mathcal{T}_k + \mathcal{T}_i + \mathcal{G}_k + \mathcal{G}_i + 4 {\mathcal{G}} \sum_{n=3}^\infty \left( \frac{1}{\delta^{3/2} \sqrt{N}} - 1  \right)^{n-3} \\
		& \lesssim \mathcal{T}_k + \mathcal{T}_i + \mathcal{G}_k + \mathcal{G}_i + 4 {\mathcal{G}} \frac{\delta^{3/2} \sqrt{N}}{2 \delta^{3/2} \sqrt{N} -1}.
	\end{align*}
	In order to obtain bounds on $U_k^n$, first consider the term
	\begin{align*}
		\sum_m \norm{[T^n]_{km} - [\bar{T}^n]_{km}},
	\end{align*}
	inductively. The case $n=1$ is simply equivalent to $\mathcal{T}_k$.  For $n=2$, we have
	\begin{align*}
		\norm{\sum_m \sum_l T_{kl}T_{lm} - \bar{T}_{kl} \bar{T}_{lm}}  & \leq \sum_m \norm{ \sum_l T_{kl} - \bar{T}_{kl} T_{lm} } + \norm{\sum_l \bar{T}_{kl} \sum_m  (T_{lm} - \bar{T}_{lm})} \\
		& = \sum_m \norm{\sum_l T_{kl} - \bar{T}_{kl} T_{lm} } \\
		& \leq \left(\frac{1}{\delta^{3/2} \sqrt{N}} - 1 \right) \mathcal{T}_k,
	\end{align*}
	where the last inequality follows from the same reasoning as Lemma \ref{lem:philip}.
	Then for $n=3$, we have
	\begin{align*}
		\norm{\sum_m \sum_l [T^2]_{kl}T_{lm} - [\bar{T}^2]_{kl} \bar{T}_{lm}} &\leq \sum_m \sum_l \norm{[T^2]_{kl} - [\bar{T}^2]_{kl}} T_{lm} \\
		& + \norm{\sum_l [\bar{T}^2]_{kl} \sum_m  (T_{lm} - \bar{T}_{lm})} \\
		& \leq \left(\frac{1}{\delta^{3/2} \sqrt{N}} - 1 \right) \sum_l \norm{[T^2]_{kl} - [\bar{T}^2]_{kl}} \\
		& \leq \left(\frac{1}{\delta^{3/2} \sqrt{N}} - 1 \right)^2 \mathcal{T}_k.
	\end{align*}
	Again, by induction for $n \geq 1$,
	\begin{align*}
		\sum_m \norm{[T^n]_{km} - [\bar{T}^n]_{km}} \leq \left(\frac{1}{\delta^{3/2} \sqrt{N}} - 1 \right)^{n-1} \mathcal{T}_k.
	\end{align*}
	Turning to $U_k^n$, for $n \geq 2$ (clearly $U_k^1 = 0$),
	\begin{align*}
		U_k^n  \leq \sum_l \norm{[T^{n-1}]_{kl} - [\bar{T}^{n-1}]_{kl}}
	\end{align*}
	so that
	\begin{align*}
		\sum_{n=1}^\infty  U_k^n &\leq \sum_{n=1}^\infty \sum_l \norm{[T^{n}]_{kl} - [\bar{T}^{n}]_{kl}}  \\
		& \leq \sum_l \norm{T_{kl} - \bar{T}_{kl}} \sum_{n=0}^\infty  \left(\frac{1}{\delta^{3/2} \sqrt{N}} - 1 \right)^{n} \\
		& \leq \mathcal{T}_k \frac{\delta^{3/2} \sqrt{N}}{2 \delta^{3/2} \sqrt{N} -1}.
	\end{align*}
	The analogous result holds for $U_i^n$.  Combining all gives
	\begin{align*}
		\sum_{n=1}^\infty \norm{S_k^n -S_i^n} & \leq  \norm{h}_\infty \left( (\mathcal{T}_k + \mathcal{T}_i + 4 {\mathcal{G}} ) \frac{\delta^{3/2} \sqrt{N}}{2 \delta^{3/2} \sqrt{N} -1}  +  \mathcal{T}_k + \mathcal{T}_i + \mathcal{G}_k + \mathcal{G}_i  \right).
	\end{align*}
	The analysis for $I_k^n$ and $I_i^n$ follows in a similar fashion, starting with $n=1$,
	\begin{align*}
		\norm{I_k^1 - I_i^1} &\leq \norm{\sum_m  \bar{T}_{km} (h(X^m) - h(\bar{X}^m)) -  \bar{T}_{im} (h(X^m) - h(\bar{X}^m)) } \\
		& \leq \frac{2L_h}{(\delta N)^{3/2}} \sum_m \norm{X^m - \bar{X}^m},
	\end{align*}
	and for $n=2$,
	\begin{align*}
		\norm{I_k^2 - I_i^2} & \leq \frac{L_h}{(\delta N)^{3/2}} \sum_m \norm{X^m - \bar{X}^m} \left( \sum_l \norm{ \bar{T}_{kl}  - \bar{\pi}_l} +\sum_l \norm{ \bar{T}_{il}  - \bar{\pi}_l} \right) \\
		& \leq \frac{4 L_h}{(\delta N)^{3/2}}  \sum_m \norm{X^m - \bar{X}^m},
	\end{align*}
	and for $n=3$,
	\begin{align*}
		\norm{I_k^3 - I_i^3} & \leq \frac{L_h}{(\delta N)^{3/2}} \sum_m \norm{X^m - \bar{X}^m} \left( \sum_l \norm{ [\bar{T}^2]_{kl}  - \bar{\pi}_l} +\sum_l \norm{ [\bar{T}^2]_{il}  - \bar{\pi}_l} \right),
	\end{align*}
	so that by induction, for $n \geq 2$
	\begin{align*}
		\norm{I_k^n - I_i^n} & \leq \frac{L_h}{(\delta N)^{3/2}} \sum_m \norm{X^m - \bar{X}^m} \left( \sum_l \norm{ [\bar{T}^{n-1}]_{kl}  - \bar{\pi}_l} +\sum_l \norm{ [\bar{T}^{n-1}]_{il}  - \bar{\pi}_l} \right) \\
		& \leq \frac{4 L_h}{(\delta N)^{3/2}}  \left( \frac{1}{\delta^{3/2} \sqrt{N}} - 1  \right)^{n-2} \sum_m \norm{X^m - \bar{X}^m},
	\end{align*}
	using (\ref{eq:boundTtostat}). Combining all, we have
	\begin{align*}
		\norm{r_k - r_i - \bar{r}_i + \bar{r}_k} &\leq \epsilon \left( 2L_h \norm{X^k - \bar{X}^k} + 2L_h \norm{X^i - \bar{X}^i}   \right) + \epsilon \sum_{n=1}^\infty \norm{S_k^n -S_i^n}   \\
		& + \epsilon \sum_{n=1}^\infty \norm{I_k^n - I_n^i}  \\
		& \lesssim  \epsilon L_h \left( \norm{X^k - \bar{X}^k} + \norm{X^i - \bar{X}^i}   \right)  \\
		& +  L_h \epsilon^2\left( \frac{  \delta^{3/2} \sqrt{N}}{2 \delta^{3/2} \sqrt{N} -1}   \right) \frac{1 }{N} \sum_m \norm{X^m - \bar{X}^m}   \\
		&   +  \epsilon \norm{h}_\infty \left(\frac{ 3 \delta^{3/2} \sqrt{N}-1}{2 \delta^{3/2} \sqrt{N} - 1} \right)  (\mathcal{T}_k + \mathcal{T}_i + 4 {\mathcal{G}} ). \\
	\end{align*}
	Making use of (\ref{eq:Tcontcond}) and
	\begin{align*}
		\mathcal{G}  & = \max_j  \left(\sum_l T_{jl} \sum_m \norm{T_{lm} - \bar{T}_{lm}} \right) \\
		& \leq \frac{1}{(\delta N)^{3/2}} \sum_l \mathcal{T}_l,
	\end{align*}
	gives the result (\ref{eq:rdiffLip}).
	\qed

	\section{Proof of lemma \ref{lem:Klip}}
	\label{app:Klip}
	
	Throughout the proof, we will make use of the following bounds
	\begin{align}
		\label{eq:Tikxixkboundnew}
		T_{ik} \norm{X^i - X^k} \leq \frac{g_{\epsilon}(X^i, X^k)}{(\delta N)^{3/2}}  \norm{X^i - X^k} & \leq \frac{2 \sqrt{d \epsilon} \exp(-1/2)}{(\delta N)^{3/2}}, \\
		\label{eq:Tikxixksquaredboundnew}
		T_{ik} \norm{X^i - X^k}^2 & \leq \frac{4 d \epsilon \exp(-1)}{(\delta N)^{3/2}}, \\
		\label{eq:Thalfikxixkboundnew}
		T_{ik}^{1/2} \norm{X^i - X^k} \leq \frac{g_{2\epsilon}(X^i, X^k)}{(\delta N)^{3/4}} \norm{X^i - X^k}  &\leq \frac{2 \sqrt{2 \epsilon d} \exp(-1/2)}{(\delta N)^{3/4}}.
	\end{align}
	We work with
	\begin{align*}
		K_\epsilon(X^i_t, \rho_t^N) = \frac{1}{2 \epsilon} \sum_{k,l} T_{ik}T_{il} (r_k - r_l) (X^k - X^i),
	\end{align*}
	since it holds that $\sum_j s_{ij} = 0$ for all $i$.
	Then using Lemma \ref{lem:philip}, we have that
	\begin{align*}
		\norm{K_\epsilon(X_t^i, \rho_t^N)} & \leq \frac{1}{2 \epsilon} \sum_k T_{ik}^{1/2} \norm{r_k - r_i}T_{ik}^{1/2}\norm{X^k - X^i} \\
		& + \frac{1}{2 \epsilon} \left(\sum_{l} T_{il}\norm{r_i - r_l} \right) \left( \sum_k T_{ik} \norm{X^k - X^i} \right)  \\
		& \leq \frac{C_\phi}{2 \epsilon} \sum_k T_{ik}^{1/2} \norm{X^k - X^i}T_{ik}^{1/2}\norm{X^k - X^i} + \frac{C_\phi}{2 \epsilon} \left( \sum_{l} T_{il}\norm{X^i - X^l} \right)^2 \\
		&  \lesssim C_\phi  d,
	\end{align*}
	since when $N$ satisfies condition \ref{ass:delt}, we have that $\frac{1}{\delta^{3/2} \sqrt{N}} < 2$.  In order to prove the second result, we start from
	\begin{align*}
		&\norm{K_\epsilon(\bar{X}^i_t, \bar{\rho}_t^N) - K_\epsilon(X_t^i, \rho_t^N)}^2 \\
		& \leq \frac{1}{2 \epsilon^2} \left(\sum_{k,l} \norm{ T_{il}T_{ik}(r_k- r_i)(X^k - X^i)- \bar{T}_{il}\bar{T}_{ik}(\bar{r}_k - \bar{r}_i)(\bar{X}^k - \bar{X}^i)} \right)^2 \\
		& + \frac{1}{2 \epsilon^2} \left( \sum_{k,l} \norm{ T_{il}T_{ik}(r_i - r_l)(X^k - X^i) - \bar{T}_{il} \bar{T}_{ik}(\bar{r}_i - \bar{r}_l)(\bar{X}^k - \bar{X}^i)} \right)^2 \\
		&=: \frac{1}{2\epsilon^2} (A_1 + A_2).
	\end{align*}
	We have that
	\begin{align*}
		A_1 & \lesssim \left( \sum_{k,l} \norm{T_{il} - \bar{T}_{il}}T_{ik}\norm{r_k - r_i} \norm{X^k - X^i} \right)^2 \\
		& + \left(\sum_{k,l} \bar{T}_{il} \norm{T_{ik}(r_k- r_i)(X^k - X^i) -\bar{T}_{ik}(\bar{r}_k - \bar{r}_i)(\bar{X}^k - \bar{X}^i) } \right)^2\\
		& =: A_{11} + A_{12},
	\end{align*}
	and likewise
	\begin{align*}
		A_2 &\lesssim \left(\sum_{k,l} T_{il} \norm{r_i - r_l} \norm{T_{ik}(X^k - X^i) - \bar{T}_{ik}(\bar{X}^k - \bar{X}^i)} \right)^2 \\
		& + \left( \sum_{k,l} \norm{\bar{T}_{ik}(\bar{X}^k - \bar{X}^i)} \norm{T_{il}(r_i - r_l) - \bar{T}_{il}(\bar{r}_i - \bar{r}_l)} \right)^2 \\
		& =: A_{21} + A_{22}.
	\end{align*}
	Using (\ref{eq:Tikxixksquaredboundnew}), Lemma \ref{lem:philip} (specifically, (\ref{eq:boundrdiff})) and Jensen's inequality we have that
	\begin{align*}
		A_{11} & \lesssim \left(\sum_{l} \norm{T_{il} - \bar{T}_{il}}  \right)^2 \sum_k T_{ik}\norm{r_k - r_i}^2 \norm{X^k - X^i}^2 \\
		& \lesssim C_\gamma^2  d \epsilon  \left( \sum_{l} \norm{T_{il} - \bar{T}_{il}}  \right)^2 \\
		& =: C_\gamma^2 d \epsilon B_i.
	\end{align*}
	
	Similarly, for $A_{12}$, we have that
	\begin{align*}
		A_{12} &  \leq 2 \left(\sum_k \norm{r_k - r_i} \norm{T_{ik}(X^k - X^i) - \bar{T}_{ik} (\bar{X}^k - \bar{X}^i)} \right)^2 \\
		& + 2 \left(\sum_k \bar{T}_{ik} \norm{\bar{X}^i - \bar{X}^k} \norm{r_k - r_i - \bar{r}_k + \bar{r}_i} \right)^2 \\
		& \leq 2 C_\gamma^2 \left(\sum_k \norm{T_{ik}(X^k - X^i) - \bar{T}_{ik} (\bar{X}^k - \bar{X}^i)} \right)^2  \\
		& + \frac{8 d \epsilon \exp(-1)}{(\delta N)^{3/2}} \sum_k \norm{r_k - r_i - \bar{r}_k + \bar{r}_i}^2 \\
		& \lesssim C_\gamma^2 \left(\sum_k \norm{T_{ik}(X^k - X^i) - \bar{T}_{ik} (\bar{X}^k - \bar{X}^i)} \right)^2 +\frac{d\epsilon}{N} \sum_k \norm{r_k - r_i - \bar{r}_k + \bar{r}_i}^2   \\
		& =: C_\gamma^2 E_i + \frac{d\epsilon}{N} D_i,
	\end{align*}
	and
	\begin{align*}
		A_{21} & \leq  C_\gamma^2 \left(\sum_k \norm{T_{ik}(X^k - X^i) - \bar{T}_{ik}(\bar{X}^k - \bar{X}^i)} \right)^2  
		= C_\gamma^2 E_i,
	\end{align*}
	and using (\ref{eq:Tikxixkboundnew}) and Jensen's inequality
	\begin{align*}
		A_{22} & \leq 16 d \epsilon \exp(-1) \left( \sum_l \norm{T_{il}(r_i - r_l) - \bar{T}_{il}(\bar{r}_i - \bar{r}_l)}  \right)^2 \\
		& \lesssim  d \epsilon \left( C_\gamma^2 \left(\sum_l \norm{T_{il} - \bar{T}_{il}} \right)^2 + \frac{1}{N}\sum_l  \norm{r_i - r_l - \bar{r}_i + \bar{r}_l}^2   \right)   \\
		&= d \epsilon \left( C_\gamma^2 B_i + \frac{1}{N} D_i \right).
	\end{align*}
	Combining all and using  Lemma \ref{lem:lipT} and \ref{lem:rrbarlip} gives the result.
	\qed

	\section{Proof of lemma \ref{ProofSLLN}}
	\label{app:ProofSLLN}
	For notational ease, we drop the $t$ subscript throughout the proof and consider $t < \zeta_\delta^N$.  Start by analysing the summand in (\ref{eq:sllntest})
	\begin{align*}
		\norm{K_\epsilon(\bar{X}^i, \bar{\rho}^N) - \bar{K}_\epsilon(\bar{X}^i, \bar{\rho})}^2  & = \frac{1}{4 \epsilon^2} \norm{  \sum_k \bar{T}_{ik} (\bar{X}^k - \hat{\bar{X}}^i) (\bar{r}^k - \hat{\bar{r}}^i)   -    \int (\mathcal{R}_\epsilon(y) - {\hat{\mathcal{R}_\epsilon}}(\bar{X}^i)) (y-\hat{y}) p_\epsilon (\bar{X}^i,y) dy} \\
		& =: \frac{1}{4 \epsilon^2} \norm{ \sum_k T_{ik}f_{ik}  -    \int f(y) p_\epsilon (\bar{X}^i,y) dy },
	\end{align*}
	where $\hat{y} := \int y p_\epsilon(\bar{X}^i,y)dy$ (likewise for $\hat{\mathcal{R}}_\epsilon$) and $\hat{X}^i := \sum_{k} T_{ik} X^k$ and likewise for $\hat{r}^i$.  For ease of notation, define $f_{ik} :=  (X^k - \hat{X}^i) (r^k - \hat{r}^i) $ and $f:= (\mathcal{R}_\epsilon - \hat{\mathcal{R}_\epsilon}) (y-\hat{y})$. Now to analyse the term in brackets, we have
	\begin{align*}
		\sum_k T_{ik}f_{ik}  -    \int f(y) p_\epsilon (\bar{X}^i,y) dy & = 	\frac{\sum_k \tilde{q}_\epsilon (\bar{X}^i , \bar{X}^k ) f_{ik}}
		{\sum_k \tilde{q}_\epsilon (\bar{X}^i , \bar{X}^k ) }  -     \int f(y) p_\epsilon (\bar{X}^i,y) dy \\
		& = 	\frac{\sum_k \tilde{q}_\epsilon (\bar{X}^i , \bar{X}^k ) (f_{ik} -      \int f(y) p_\epsilon (\bar{X}^i,y) dy )}
		{\sum_k \tilde{q}_\epsilon (\bar{X}^i , \bar{X}^k ) }  \\
		& = 	\frac{\frac{1}{N}\sum_k \tilde{q}_\epsilon (\bar{X}^i , \bar{X}^k ) \int b({\bf X}, y)  p_\epsilon (\bar{X}^i,y) dy }
		{\frac{1}{N} s^i },
	\end{align*}
	where
	\begin{align*}
		b({\bf X}, y) &:= f_{ik} -  f(y) \\
		& =  (X^k - \hat{X}^i) (r^k - \hat{r}^i)  - (\mathcal{R}_\epsilon(y) - \hat{\mathcal{R}_\epsilon}(\bar{X}^i)) (y-\hat{y}).
	\end{align*}
	From now on, we consider the expectation and $t < \zeta_\delta^N$.  Using the fact that $s^i \geq \delta N$ it holds that
	\begingroup
	\allowdisplaybreaks
	\begin{align*}
		&	\mathbb{E} \left[\norm{\frac{\frac{1}{N}\sum_k \tilde{q}_\epsilon (\bar{X}^i , \bar{X}^k ) \int b({\bf X}, y)  p_\epsilon (\bar{X}^i,y) dy }
			{\frac{1}{N} s^i } }^2 \right] \\
		& \leq  \frac{1}{\delta^2}  \mathbb{E} \left[\norm{\frac{1}{N}\sum_k \tilde{q}_\epsilon (\bar{X}^i , \bar{X}^k ) \int b({\bf X}, y)  p_\epsilon (\bar{X}^i,y) dy }^2\right] \\
		& =  \frac{1}{\delta^2}  \mathbb{E} \left[\norm{\frac{1}{N}\sum_k  \frac{g_\epsilon(\bar{X}^i, \bar{X}^k)}{(\sum_l g_\epsilon(\bar{X}^k, \bar{X}^l))^{1/2}} \int b({\bf X}, y)  p_\epsilon (\bar{X}^i,y) dy }^2\right]  \\
		& =  \frac{1}{\delta^2N}  \mathbb{E} \left[\norm{\frac{1}{N}\sum_k  \frac{g_\epsilon(\bar{X}^i, \bar{X}^k)}{( \bar{z}_\epsilon^k)^{1/2}} \int b({\bf X}, y)  p_\epsilon (\bar{X}^i,y) dy }^2\right] \\
		& \leq \frac{2}{\delta^2 N} \mathbb{E} \left[\norm{\frac{1}{N}\sum_k  \frac{g_\epsilon(\bar{X}^i, \bar{X}^k)}{(g_\epsilon * \rho(\bar{X}^k))^{1/2}} \int b({\bf X}, y)  p_\epsilon (\bar{X}^i,y) dy}^2\right]
		+ \frac{2}{\delta^2 N}  \mathbb{E}\\& \left[\norm{\frac{1}{N}\sum_k  \int b({\bf X}, y)  p_\epsilon (\bar{X}^i,y) dy  g_\epsilon(\bar{X}^i, \bar{X}^k) \left( \frac{1}{(\frac{1}{N} \bar{z}_\epsilon^k)^{1/2}}-  \frac{1}{(g_\epsilon * \rho(\bar{X}^k))^{1/2}}     \right) }^2\right]  \\
		& =: \frac{2}{\delta^2 N} (I_1 + I_2),
	\end{align*}
	\endgroup
	where $\bar{z}_\epsilon^k := \frac{1}{N} \sum_l g_\epsilon(\bar{X}^k, \bar{X}^l)$.  First note that $\int b({\bf X}, y)  p_\epsilon (\bar{X}^i,y) dy  g_\epsilon(\bar{X}^i, \bar{X}^k) $ is uniformly bounded for all $i,k$ since
	\begin{align*}
		\norm{g_\epsilon(\bar{X}^i, \bar{X}^k) \int b({\bf X}, y)  p_\epsilon (\bar{X}^i,y) dy  } \leq \norm{f_{ik}} g_\epsilon(\bar{X}^i, \bar{X}^k) + \norm{ \int f(y) p_\epsilon(\bar{X}^i, y) dy},
	\end{align*}
	and using Lemma \ref{lem:philip}, (\ref{eq:Tbounds}) and (\ref{eq:Tikxixkbound}), we have that for all $i,k$,
	\begin{align*}
		\norm{f_{ik} }g_\epsilon(\bar{X}^i, \bar{X}^k) \lesssim C_\gamma d (\epsilon \vee \sqrt{\epsilon}).
	\end{align*}
	Secondly, under assumption (\ref{ass:existsoln}) and  (\ref{ass:nubound}) with $h$ uniformly bounded, we have
	\begin{align*}
		&\norm{ \int f(y) p_\epsilon(\bar{X}^i, y) dy} \\
		&= \norm{\int \left(\mathcal{R}_\epsilon(y) - \int \mathcal{R}_\epsilon (y) p_\epsilon(\bar{X}^i,y) dy \right)  \left(y - \int y p_\epsilon(\bar{X}^i,y) dy \right) p_\epsilon(\bar{X}^i,y) dy } \\
		& \leq \int  \norm{\mathcal{R}_\epsilon(y) - \int \mathcal{R}_\epsilon (y) p_\epsilon(\bar{X}^i,y) dy}  \norm{y - \int y p_\epsilon(\bar{X}^i,y) dy} p_\epsilon(\bar{X}^i,y) dy \\
		& \lesssim (\norm{\phi_\epsilon}_\infty + \epsilon \norm{h}_\infty) \delta_v,
	\end{align*}
	where the boundedness of $\phi_\epsilon$ follows from the fact that the rhs of the fixed point equation (\ref{eq:phidiffmap}) is uniformly bounded.  Combining gives
	\begin{align*}
		\norm{g_\epsilon(\bar{X}^i, \bar{X}^k) \int b({\bf X}, y)  p_\epsilon (\bar{X}^i,y) dy  } & \leq  (\norm{\phi_\epsilon}_\infty + \epsilon \norm{h}_\infty) \delta_v + C_\gamma d (\epsilon \vee \sqrt{\epsilon}) \\
		& =: D_b.
	\end{align*}
	Starting with $I_1$, we have
	\begin{align*}
		I_1	& \leq   D_b^2 \mathbb{E} \left[\norm{\frac{1}{N}\sum_k  \frac{1}{(g_\epsilon * \rho(\bar{X}^k))^{1/2}} }^2\right] \\
		& \leq D_b^2 \frac{1}{N}\sum_k\mathbb{E} \left[ \frac{1}{g_\epsilon * \rho(\bar{X}^k)}\right] \\
		& \leq \frac{D_b^2}{\delta},
	\end{align*}
	since for $t < \zeta_\delta^N$, it holds that for all $k$
	\begin{align}
		\label{eq:gconvbound}
		g_\epsilon \ast \rho(\bar{X}^k) \geq \delta.
	\end{align}
	Then for $I_2$,
	\begin{align*}
		I_2 & \leq  D_b^2 \mathbb{E} \left[ \frac 1{N}
		\sum_k \frac{\left| \bar{z}_\epsilon^k  - g_\epsilon \ast \rho (\bar{X}^k)\right|}{\bar{z}_\epsilon^k
			g_\epsilon\ast\rho    (\bar{X}^k ) } \right] \\
		& \leq \frac{ D_b^2}{\delta} \mathbb{E} \left[ \frac 1{N}
		\sum_k \frac{\left| \bar{z}_\epsilon^k  - g_\epsilon \ast \rho (\bar{X}^k)\right|}{
			g_\epsilon\ast\rho    (\bar{X}^k ) } \right] \\
		& \leq \frac{ D_b^2}{\delta} \frac{1}{N} \sum_k  \mathbb{E} \left[ \frac 1{g_\epsilon\ast\rho (\bar{X}^k )^2}\right]^{\frac 12}  \mathbb{E} \left[
		\norm{\bar{z}_\epsilon^k - g_\epsilon \ast \rho (\bar{X}^k)}^2 \right]^{\frac 12} \\
		& \leq \frac{ D_b^2}{\delta} \frac{1}{N} \sum_k  \mathbb{E} \left[ \frac 1{g_\epsilon\ast\rho(\bar{X}^k )^2}\right]^{\frac 12}  \left(
		\frac{1}{N} \sum_l \mathbb{E} \left[ \norm{ g_\epsilon(\bar{X}^k, \bar{X}^l) - g_\epsilon \ast \rho (\bar{X}^k) }^2 \right] \right)^{\frac 12} \\
		& \lesssim \frac{ D_b^2}{\delta} \frac{1}{N} \sum_k  \mathbb{E} \left[ \frac 1{g_\epsilon\ast\rho (\bar{X}^k )^2}\right]^{\frac 12},
	\end{align*}
	where the last inequality arises from the fact that for all $l,k$, it holds that
	\begin{align*}
		\mathbb{E} \left[ \norm{ g_\epsilon(\bar{X}^k, \bar{X}^l) - g_\epsilon \ast \rho (\bar{X}^k) }^2 \right] & \leq 2 \mathbb{E} \left[ \norm{ g_\epsilon(\bar{X}^k, \bar{X}^l)}^2 \right] + 2 \mathbb{E} \left[  \norm{g_\epsilon \ast \rho (\bar{X}^k) }^2 \right] \\
		& \leq 4.
	\end{align*}
	To analyse the remaining expectation, we have that using (\ref{eq:gconvbound}) for every $k$,
	\begin{align*}
		\mathbb{E} \left[ \frac 1{g_\epsilon\ast\rho (\bar{X}^k )^2}\right]  \leq \frac{1}{\delta^2}.
	\end{align*}
	Combining leads to
	\begin{align*}
		I_2 \leq \frac{D_b^2}{\delta^2}.
	\end{align*}
	the result then follows easily.
	\qed

	\section{Some useful theorems}
	
	\begin{lemma}
		\label{lem:covPI}
		\textbf{Generalised Brascamp-Lieb inequality, theorem 1.1 in \cite{Carlen2013}} Let $g,h: \mathbb{R}^d \rightarrow \mathbb{R}$ be square integrable locally Lipschitz functions and suppose $d\mu = \exp(-V)dx$ is a strictly log-concave measure. It then holds that
		\begin{align*}
			\norm{Cov_\mu (g,h)} & \leq \norm{Hess(V)^{-1/2} \nabla g}_2 \norm{Hess(V)^{-1/2} \nabla h} \\
			& = \left( \int \nabla g^T Hess(V)^{-1} \nabla g d\mu  \right)^{1/2} \left( \int \nabla h^T Hess(V)^{-1} \nabla h d\mu  \right)^{1/2},
		\end{align*}	
		where the norms on the rhs refer to the euclidean norm (although the result holds for more general $L^p$ norms.)
		
		Suppose furthermore that $V$ is uniformly convex, i.e. $\inf_x y^T Hess(V)(x) y \geq \kappa |y|^2$ for all $y \in \mathbb{R}^d$, then  	
		\begin{align*}
			\norm{Cov_\mu (g,h)} & \leq \frac{1}{\kappa} \left( \int \norm{\nabla g}^2 d\mu  \right)^{1/2}  \left( \int \norm{\nabla h}^2 d\mu  \right)^{1/2}.
		\end{align*}
	\end{lemma}

\bibliographystyle{halpha}
\bibliography{literature_DMFPF.bib}  

\begin{thebibliography}{LMMR15}
\expandafter\ifx\csname url\endcsname\relax
  \def\url#1{\texttt{#1}}\fi
\expandafter\ifx\csname doi\endcsname\relax
  \def\doi#1{\burlalt{doi:#1}{http://dx.doi.org/#1}}\fi
\expandafter\ifx\csname urlprefix\endcsname\relax\def\urlprefix{URL }\fi
\expandafter\ifx\csname href\endcsname\relax
  \def\href#1#2{#2}\fi
\expandafter\ifx\csname burlalt\endcsname\relax
  \def\burlalt#1#2{\href{#2}{#1}}\fi

\bibitem[BdP18]{Bishop_P_2018}
A.~N. Bishop, P.~{del Moral}, and S.~Pathiraja.
\newblock {Perturbations and projections of Kalman-Bucy semigroups}.
\newblock {\em Stochastic Processes and their Applications}, 128(9), 2018,
  \burlalt{1701.05978}{http://arxiv.org/abs/1701.05978}.
\newblock \doi{10.1016/j.spa.2017.10.006}.

\bibitem[BG16]{Berntorp2016}
K.~Berntorp and P.~Grover.
\newblock {Data-driven gain computation in the feedback particle filter}.
\newblock In {\em Proceedings of the American Control Conference}, pages
  2711--2716, jul 2016.
\newblock \doi{10.1109/ACC.2016.7525328}.

\bibitem[BGM10]{Bolley2010}
F.~Bolley, A.~Guillin, and F.~Malrieu.
\newblock {Trend to equilibrium and particle approximation for a weakly
  selfconsistent Vlasov-Fokker-Planck equation}.
\newblock {\em ESAIM: Mathematical Modelling and Numerical Analysis},
  44(5):867--884, 2010, \burlalt{0906.1417}{http://arxiv.org/abs/0906.1417}.
\newblock \doi{10.1051/m2an/2010045}.

\bibitem[BH16]{Berry2016}
T.~Berry and J.~Harlim.
\newblock {Variable bandwidth diffusion kernels}.
\newblock {\em Applied and Computational Harmonic Analysis}, 40(1):68--96,
  2016, \burlalt{1406.5064}{http://arxiv.org/abs/1406.5064}.
\newblock \doi{10.1016/j.acha.2015.01.001}.

\bibitem[BL76]{Brascamp1976}
H.~J. Brascamp and E.~H. Lieb.
\newblock {On extensions of the Brunn-Minkowski and Pr{\'{e}}kopa-Leindler
  theorems, including inequalities for log concave functions, and with an
  application to the diffusion equation}.
\newblock {\em Journal of Functional Analysis}, 22(4):366--389, 1976.
\newblock \doi{10.1016/0022-1236(76)90004-5}.

\bibitem[BR12]{Bergemann2012}
K.~Bergemann and S.~Reich.
\newblock {An ensemble Kalman-Bucy filter for continuous data assimilation}.
\newblock {\em Meteorologische Zeitschrift}, 21(3):213--219, 2012.
\newblock \doi{10.1127/0941-2948/2012/0307}.

\bibitem[CCeL13]{Carlen2013}
E.~A. Carlen, D.~Cordero-erausquin, and E.~H. Lieb.
\newblock {Asymmetric covariance estimates of Brascamp-Lieb type and related
  inequalities for log-concave measures}.
\newblock {\em Annales de l'Institut Henri Poincare}, 49(1):1--12, 2013.
\newblock \doi{10.1214/11-AIHP462}.

\bibitem[CL06]{Coifman2006}
R.~R. Coifman and S.~Lafon.
\newblock {Diffusion maps}.
\newblock {\em Applied Computational Harmonic Analysis}, 21:5--30, 2006.
\newblock \doi{10.1016/j.acha.2006.04.006}.

\bibitem[CX10]{sr:crisan10}
D.~Crisan and J.~Xiong.
\newblock {Approximate McKean-Vlasov representations for a class of SPDEs}.
\newblock {\em Stochastics}, 82, 2010.
\newblock \doi{10.1080/17442500902723575}.

\bibitem[dKT17]{DelMoral2017}
P.~{del Moral}, A.~Kurtzmann, and J.~Tugaut.
\newblock {On the Stability and the Uniform Propagation of Chaos of a Class of
  Extended Ensemble Kalman--Bucy Filters}.
\newblock {\em SIAM Journal on Control and Optimization}, 55(1):119--155, 2017.
\newblock \doi{10.1137/16M1087497}.

\bibitem[dWRS18]{deWiljesStannat2018}
J.~de~Wiljes, S.~Reich, and W.~Stannat.
\newblock {Long-time stability and accuracy of the ensemble Kalman-Bucy filter
  for fully observed processes and small measurement noise}.
\newblock {\em SIAM J. Applied Dynamical Systems}, 17(2):1152--1181, 2018,
  \burlalt{1612.06065}{http://arxiv.org/abs/1612.06065}.
\newblock \doi{10.1137/17M1119056}.

\bibitem[dWT20]{deWiljes2020}
J.~de~Wiljes and X.~T. Tong.
\newblock {Analysis of a localised nonlinear ensemble Kalman Bucy filter with
  complete and accurate observations}.
\newblock {\em Nonlinearity}, 33(9):4752--4782, 2020.
\newblock \doi{10.1088/1361-6544/ab8d14}.

\bibitem[Eve03]{sr:evensen03}
G.~Evensen.
\newblock {The Ensemble Kalman Filter: theoretical formulation and practical
  implementation}.
\newblock {\em Ocean Dynamics}, 53(4):343--367, 2003.
\newblock \doi{10.1007/s10236-003-0036-9}.

\bibitem[Eve06]{sr:evensen}
G.~Evensen.
\newblock {\em Data assimilation. {T}he ensemble Kalman filter}.
\newblock Springer-Verlag, New York, 2006.

\bibitem[EvLP00]{Evensen2000}
G.~Evensen and van Leeuwen~P.J.
\newblock {An Ensemble Kalman Smoother for Nonlinear Dynamics}.
\newblock {\em Monthly Weather Review}, 128:1852--1867, 2000.

\bibitem[LMMR15]{Laugesen2015}
R.~S. Laugesen, P.~G. Mehta, S.~P. Meyn, and M.~Raginsky.
\newblock {Poissons Equation in Nonlinear Filtering}.
\newblock {\em SIAM Journal on Control and Optimization}, 53(1):501--525, 2015.
\newblock \doi{10.1137/13094743X}.

\bibitem[LMT11]{LeGland2009}
F.~{Le Gland}, V.~Monbet, and V.~Tran.
\newblock {Large sample asymptotics for the ensemble Kalman filter}.
\newblock {\em Oxford Handbook of Nonlinear Filtering}, 2011.

\bibitem[LS20]{Lange2020}
T.~Lange and W.~Stannat.
\newblock {Mean field limit of Ensemble Square Root Filters - Discrete and
  continuous time}.
\newblock {\em Foundations of Data Science}, pages 1--26, 2020,
  \burlalt{2011.10516}{http://arxiv.org/abs/2011.10516}.
\newblock \doi{10.3934/fods.2021003}.

\bibitem[LS21]{Lange2019}
T.~Lange and W.~Stannat.
\newblock {On the continuous time limit of the ensemble kalman filter}.
\newblock {\em Mathematics of Computation}, 90:233--265, 2021,
  \burlalt{1901.05204}{http://arxiv.org/abs/1901.05204}.
\newblock \doi{10.1090/mcom/3588}.

\bibitem[MT16]{Majda2016}
A.J. Majda and X.T. Tong.
\newblock {Robustness and Accuracy of finite Ensemble Kalman filters in large
  dimensions}.
\newblock {\em arXiv}, page~40, jun 2016,
  \burlalt{1606.0932}{http://arxiv.org/abs/1606.0932}.
\newblock \urlprefix\url{http://arxiv.org/abs/1606.0932}.

\bibitem[OTM20]{Olmez2020}
S.~Y. Olmez, A.~Taghvaei, and P.G. Mehta.
\newblock {Deep FPF : Gain function approximation in high-dimensional setting}.
\newblock 2020,
  \burlalt{arXiv:2010.01183v1}{http://arxiv.org/abs/arXiv:2010.01183v1}.

\bibitem[Pat20]{P_wongzakaiarxiv2020}
S.~Pathiraja.
\newblock {L2 convergence of smooth approximations of Stochastic Differential
  Equations with unbounded coefficients}.
\newblock {\em arXiv arXiv:2011.13009}, 2020.
\newblock \urlprefix\url{https://arxiv.org/abs/2011.13009}.

\bibitem[PRS21]{Pathiraja2021}
S.~Pathiraja, S.~Reich, and W.~Stannat.
\newblock {McKean-Vlasov SDEs in nonlinear filtering}.
\newblock {\em SIAM J. Control and Optimization}, 59(6):4188--4215, 2021,
  \burlalt{2007.12658}{http://arxiv.org/abs/2007.12658}.
\newblock \doi{10.1137/20m1355197}.

\bibitem[RDM16]{Radhakrishnan2016}
A.~Radhakrishnan, A.~Devraj, and S.~Meyn.
\newblock {Learning techniques for feedback particle filter design}.
\newblock In {\em 2016 IEEE 55th Conference on Decision and Control (CDC)},
  pages 5453--5459, dec 2016.
\newblock \doi{10.1109/CDC.2016.7799106}.

\bibitem[Rei10]{sr:reich10}
S.~Reich.
\newblock {A dynamical systems framework for intermittent data assimilation}.
\newblock {\em BIT Numerical Mathematics}, 51:235--249, 2010.
\newblock \doi{10.1007/s10543-010-0302-4}.

\bibitem[SS16]{Schillings2016}
C.~Schillings and A.M. Stuart.
\newblock {Analysis of the ensemble kalman filter for inverse problems}.
\newblock {\em arxiv}, jan 2016,
  \burlalt{1602.0202}{http://arxiv.org/abs/1602.0202}.
\newblock \urlprefix\url{http://arxiv.org/abs/1602.0202}.

\bibitem[TMK16]{TongKelly2016}
X.~T. Tong, A.~J. Majda, and D.~Kelly.
\newblock {Nonlinear stability and ergodicity of ensemble based Kalman
  filters}.
\newblock {\em Nonlinearity}, 29(2):657--691, feb 2016.
\newblock \doi{10.1088/0951-7715/29/2/657}.

\bibitem[TMM20]{Taghvaei2019}
A.~Taghvaei, P.~G. Mehta, and S.~P. Meyn.
\newblock {Diffusion map-based algorithm for gain function approximation in the
  feedback particle filter}.
\newblock {\em SIAM-ASA Journal on Uncertainty Quantification},
  8(3):1090--1117, 2020, \burlalt{1902.07263}{http://arxiv.org/abs/1902.07263}.
\newblock \doi{10.1137/19M124513X}.

\bibitem[WR21]{Wormell2021}
C.~L. Wormell and S.~Reich.
\newblock {Spectral Convergence of Diffusion Maps: Improved Error Bounds and an
  Alternative Normalization}.
\newblock {\em SIAM Journal on Numerical Analysis}, 59(3):1687--1734, 2021,
  \burlalt{2006.02037}{http://arxiv.org/abs/2006.02037}.
\newblock \doi{10.1137/20m1344093}.

\bibitem[YMM11]{Yangmfc2011}
T.~Yang, P.~G. Mehta, and S.~P. Meyn.
\newblock {A mean-field control-oriented approach to particle filtering}.
\newblock {\em Proceedings of the American Control Conference}, pages
  2037--2043, 2011.
\newblock \doi{10.1109/acc.2011.5991422}.

\bibitem[YMM13]{Yang2013}
T.~Yang, P.~G. Mehta, and S.~P. Meyn.
\newblock {Feedback Particle Filter}.
\newblock {\em IEEE Transactions on Automatic Control}, 58(10):2465--2480,
  2013.
\newblock \doi{10.1109/TAC.2013.2258825}.

\end{thebibliography}

\end{document}